\documentclass{article}

\usepackage[utf8]{inputenc}
\usepackage[T1]{fontenc}
\usepackage{hyperref}
\usepackage{url}
\usepackage{booktabs}
\usepackage{amsmath,amssymb,amsthm,amsfonts}
\usepackage{mathtools}
\usepackage{mathrsfs}
\usepackage{microtype}
\usepackage{xcolor}
\usepackage{graphicx}
\usepackage{subcaption}
\usepackage{nicefrac}

\hypersetup{
  colorlinks=true,
  linkcolor=black!60!blue,
  citecolor=black!60!blue,
  urlcolor=black!60!blue,
}

\theoremstyle{plain}
\newtheorem{theorem}{Theorem}

\newtheorem{conjecture}[theorem]{Conjecture}
\theoremstyle{definition}

\theoremstyle{remark}
\newtheorem{remark}[theorem]{Remark}

\newcommand{\R}{\mathbb R}

\newcommand{\N}{\mathbb N}

\newcommand{\F}{\mathcal F}
\newcommand{\sech}{\operatorname{sech}}

\newcommand{\be}{\begin{equation}}
\newcommand{\ee}{\end{equation}}

\title{Neural Discovery of Strichartz Extremizers}

\author{%
  Nicol\'as Valenzuela \\
  Departamento de Ingenier\'ia Matem\'atica\\
  Universidad de Chile \\
  Santiago, Chile \\
  \texttt{nvalenzuela@dim.uchile.cl}
  \\
  \\
  Ricardo Freire \\
  Departamento de Ingenier\'ia Matem\'atica\\
  Universidad de Chile \\
  Santiago, Chile \\
  \texttt{rfreire@dim.uchile.cl}
  \\
  \\
  Claudio Mu\~noz \\
  DIM \& CMM (UMI 2807 CNRS)\\
  Universidad de Chile \\
  Santiago, Chile \\
  \texttt{cmunoz@dim.uchile.cl}
}

\begin{document}

\maketitle

\begin{abstract}
Strichartz inequalities are a cornerstone of the modern theory of dispersive PDEs, but their extremizers are known explicitly only in a handful of sharp cases. The non-convexity of the underlying functional makes the problem hard, and to our knowledge no systematic numerical attack has been attempted. We propose a simple neural-network-based pipeline that searches for extremizers as critical points of the Strichartz ratio, and apply it in three settings. First, on the Schr\"odinger group we recover the Gaussian extremizers of Foschi and Hundertmark--Zharnitsky in dimensions $d=1,2$ to within $10^{-3}$ relative error, with no analytical prior. Second, on $59$ further admissible pairs in $d=1$ where the answer is conjectural, the method consistently finds Gaussians, supporting the conjecture that Gaussians are the universal extremizers in the admissible range. Third, on the critical Airy--Strichartz inequality at $\gamma=1/q$, where existence is open, the optimization does not converge to any $L^2$ profile: instead, the iterates organize themselves as mKdV breathers $B(0,\cdot;\alpha,1,0,0)$ with growing internal frequency $\alpha$, and the discovered ratio approaches the Frank--Sabin universal lower bound $\widetilde A_{q,r}$ from below with a power-law gap $\sim\alpha^{-0.9}$. We confirm the same picture with an independent Hermite-basis ansatz. We propose a precise conjecture: the supremum equals $\widetilde A_{q,r}$ and is approached, but not attained, along the breather family. The pipeline thus serves both as a validator on known cases and as a discovery tool when no extremizer exists.
\end{abstract}

\section{Introduction}
\label{sec:intro}

\paragraph{Strichartz estimates.} Strichartz estimates constitute the analytic backbone of the theory of dispersive partial differential equations. Originally introduced by \cite{Strichartz77}, they provide quantitative space-time integrability bounds for solutions of linear dispersive flows, turning purely dispersive decay into an integrable statement in mixed Lebesgue spaces $L^q_tL^r_x$. For the free Schr\"odinger group on $\R^d$, the prototypical estimate reads
\be\label{Str-Sch}
\|e^{it\Delta}u_0\|_{L^{2+4/d}(\R^{d+1})} \;\le\; S_d\,\|u_0\|_{L^2(\R^d)},
\ee
and analogous inequalities hold for the wave, Klein--Gordon, and Airy groups. Such estimates, together with their refinements and bilinear relatives \cite{LinaresPonce,KPV0,KPV1}, are the indispensable tool behind essentially all modern well-posedness theory for nonlinear dispersive models.

The question of \emph{optimal constants} and the \emph{existence of extremizers} in these inequalities, however, is a much more delicate matter. Even in the simplest Schr\"odinger setting the answer is non-trivial: \cite{Foschi07} proved that extremizers of \eqref{Str-Sch} exist in dimensions $d=1$, $p=6$ and $d=2$, $p=4$, and that they are complex Gaussians of the form $A\,e^{-bx^2+cx+d}$. This fact was revisited by \cite{HundertmarkZharnitsky06} through a Rayleigh quotient identity, by \cite{BBC H09} using heat-flow monotonicity, and it fits into a broader line of work on sharp Fourier restriction \cite{ChristQuilodran14,ChristShao12,GoncalvesNegro22}. Outside of this sharp window, however, the situation changes drastically: when the exponent is non-sharp, or outside a certain favorable setting, extremizers typically do not exist or are unknown to exist. Maximizing sequences then develop concentration and escape to infinity via the symmetry group of the equation, a phenomenon quantified by profile decomposition arguments \cite{Shao09,Kunze03}.

\paragraph{The Airy--Strichartz inequality.} The focus of this paper is the Airy group $e^{-t\partial_x^3}$, linearization of the generalized Korteweg--de Vries (gKdV) equations $\partial_t u+\partial_x(\partial_x^2 u+u^p)=0$ at zero. Following \cite{KPV0,KPV1}, for admissible triples $(q,r,\gamma)$ with $2\le q,r<\infty$, $-\gamma+3/q+1/r=1/2$, and $-1/2<\gamma\le 1/q$, one has
\be\label{AS}
\|\,|D_x|^\gamma e^{-t\partial_x^3}u_0\|_{L^q_tL^r_x(\R\times\R)} \;\le\; A_{q,r}\,\|u_0\|_{L^2_x(\R)},
\ee
where $|D_x|^\gamma = \F^{-1}(|\xi|^\gamma\F(\cdot))$. The sharp constant is
\be\label{Aqr}
A_{q,r} \;:=\; \sup_{u_0\in L^2_x(\R)\setminus\{0\}} \frac{\|\,|D_x|^\gamma e^{-t\partial_x^3}u_0\|_{L^q_tL^r_x}}{\|u_0\|_{L^2_x}}.
\ee
The critical case of \eqref{AS}, and the one we study here, is $\gamma=1/q$. In this regime, \cite{FS18} proved a striking precompactness criterion: maximizing sequences for $A_{q,r}$ are precompact modulo the symmetry group of the Airy equation \emph{if and only if}
\be\label{FS-bound}
A_{q,r} \;>\; \widetilde A_{q,r} \;:=\; 3^{-1/q}\,a_r\,S_{1,q,r},
\qquad
a_r=\frac{2^{1/2}}{\pi^{1/(2r)}}\!\left(\frac{\Gamma(\tfrac{r+1}{2})}{\Gamma(\tfrac{r+2}{2})}\right)^{1/r}>1,
\ee
where $S_{1,q,r}$ denotes the sharp Schr\"odinger--Strichartz constant. In particular, $\widetilde A_{q,r}$ is a universal lower bound for $A_{q,r}$ in terms of the Schr\"odinger one, but the existence of an actual $L^2$-attained extremizer is only guaranteed by strict inequality \eqref{FS-bound}. The borderline regime $A_{q,r}=\widetilde A_{q,r}$ is exactly the case where an extremizer fails to exist, and the inequality is saturated only as the limit of a non-compact maximizing sequence.

\paragraph{Whether extremizers are attained is open.} Before this work, to the best of our knowledge, no conjectural value of $A_{q,r}$ and no candidate profile for a limiting maximizing sequence had been proposed in the literature at the critical exponent $\gamma=1/q$ for the Airy group. Our main contribution is precisely to propose such a candidate, informed by an AI-based discovery procedure, and to provide consistent numerical evidence that the critical Airy--Strichartz inequality lies at the borderline $A_{q,r}=\widetilde A_{q,r}$, with maximizing sequences given by mKdV breathers whose internal frequency parameter diverges.

\paragraph{Why a learning-based approach?} A natural alternative is to project the variational problem onto a fixed basis (Hermite, Fourier) and run gradient descent over the coefficients. We tried this with a $21$-term Hermite expansion (Appendix~\ref{app:hermite}) as a sanity check, and it works in the qualitative sense, but the truncation $N$ and the spatial scale of the basis bias the answer in a way that is hard to disentangle. A neural-network ansatz removes both choices: the same architecture handles a Gaussian Schr\"odinger extremizer and the highly oscillatory breather family appearing for Airy. Equally important, the training trajectory itself is informative. When an $L^2$ extremizer exists the iterates stabilise at a fixed profile; when none exists, they keep drifting along a non-compact direction in function space. This second behaviour is precisely the qualitative signature we exploit in the Airy case.

\paragraph{Contributions.} The paper has four contributions.
(i) A neural-network discovery pipeline for Strichartz-type variational problems (Section~\ref{sec:method}). The pipeline rests on an FFT-based linear propagator, a quadrature approximation of the target mixed norm, and previous PINN approximation-error theorems of ours on dispersive models~\cite{MV,ACFMV,FMV,MMP}.
(ii) Validation on the linear Schr\"odinger model. Foschi's Gaussian extremizers at the sharp exponents $(d,p)\in\{(1,6),(2,4)\}$ are recovered to within $10^{-3}$ relative error; in the non-sharp admissible pairs we tested, the discovered profile is again a Gaussian.
(iii) For the critical Airy--Strichartz inequality with $\gamma=1/q$, training does \emph{not} stabilise at an $L^2$ profile. The iterates organise themselves as mKdV breathers with growing internal frequency, and the discovered ratio approaches the Frank--Sabin lower bound $\widetilde A_{q,r}$ monotonically from below.
(iv) We collect this evidence into Conjecture~\ref{conj:main}: in the critical Airy--Strichartz inequality the supremum is not attained, equals $\widetilde A_{q,r}$, and is realized as $\lim_{\alpha\to\infty} R[B(0,\cdot;\alpha,1,0,0)]$.

\paragraph{Related work.} The search for Strichartz extremizers has a long history starting with \cite{Kunze03}, \cite{Foschi07} and \cite{HundertmarkZharnitsky06}, and developed by \cite{Carneiro09}, \cite{ChristQuilodran14,ChristShao12}, \cite{Jeavons15} and \cite{GoncalvesNegro22}, among many others; we refer to \cite{BezRogers13} for the sharp Strichartz estimate for the wave equation in the energy space, to \cite{FoschiOliveira17} for a comprehensive survey of sharp Fourier restriction theory, and to the monographs \cite{LinaresPonce,Cazenave,Tao_book} for the underlying linear and nonlinear dispersive theory. On the ML side, PINNs were introduced by \cite{Raissi,RK,RPK} and have been rigorously analysed in a dispersive PDE setting in our previous work \cite{MV,ACFMV,FMV,MMP,CMV,MR}, alongside parallel developments for radiative transfer, parabolic and elliptic equations \cite{MM22,BKM22,OMO,CSHF,JGS1}; broader operator-learning frameworks such as DeepONet \cite{LuJinKar21} and the Fourier neural operator \cite{LiKovachki21}, together with the survey by \cite{Karniadakis21}, provide a wider methodological context. These references inform the present pipeline, in particular the use of dispersive-norm losses for AI-based approximators. The application of PINNs to Strichartz \emph{variational} problems, rather than as PDE solvers, appears to be new.


\section{Preliminaries}
\label{sec:prelim}

\subsection{Strichartz estimates for dispersive groups}

Given a linear dispersive propagator $V:L^2_x(\R^d)\to C(\R;L^2_x(\R^d))$ and a Banach space $X$ of space-time functions, the sharp Strichartz constant is
\be\label{sharp}
C \;:=\; \sup_{u_0\in L^2_x,\ \|u_0\|_{L^2_x}=1} \|V(t)u_0\|_X.
\ee
The search for extremizers amounts to the minimization of $-\|V(t)u_0\|_X$ under the constraint $\|u_0\|_{L^2_x}=1$. A priori, such a minimizer may fail to exist for one of two reasons: \emph{escape to infinity}, in which a maximizing sequence is not bounded in $L^2$, or \emph{concentration}, in which the sequence loses compactness along the symmetry group of $V$. Both phenomena are well documented in the analytic literature \cite{Kunze03,Shao09,FS18}.

For the Schr\"odinger group $V(t)=e^{it\Delta}$ in $\R^d$, the sharp exponents $p=2+4/d$ correspond to $d=1$, $p=6$ and $d=2$, $p=4$, and the extremizer has the explicit Gaussian form \cite{Foschi07}. For the Airy group $V(t)=e^{-t\partial_x^3}$ in $\R$, the Strichartz family is given by \eqref{AS}, and the critical case $\gamma=1/q$ is the one in which the Frank--Sabin criterion \eqref{FS-bound} applies. 

\subsection{Admissibility and the Frank--Sabin bound}

In the critical case $\gamma=1/q$, the admissibility conditions reduce to $(q,r)$ Schr\"odinger-admissible \emph{in one dimension}, i.e., $\tfrac{1}{q}+\tfrac{1}{2r}=\tfrac14$ with $2\le q,r<\infty$. The four representative pairs we consider are $(q,r)\in\{(5,10),(6,6),(8,4),(12,3)\}$. The Frank--Sabin universal lower bound \eqref{FS-bound} produces a four-value baseline
\[
\widetilde A_{5,10}\approx 0.7926,\qquad \widetilde A_{6,6}\approx 0.7886,\qquad \widetilde A_{8,4}\approx 0.8112,\qquad \widetilde A_{12,3}\approx 0.8556.
\]
By the theorem of \cite{FS18}, extremizers for \eqref{AS} exist in these pairs if and only if the corresponding $A_{q,r}$ is \emph{strictly} larger than its respective $\widetilde A_{q,r}$. Our numerical evidence, presented in Section~\ref{sec:airy}, is consistent with equality rather than strict inequality, i.e., with \emph{non-existence of extremizers}, except in some particular cases.

\subsection{mKdV breathers}

The modified Korteweg--de Vries equation $\partial_t u+\partial_x(\partial_x^2 u+u^3)=0$ admits, for each $\alpha,\beta>0$, the explicit \emph{breather} solutions \cite{Wadati,Lamb,AM}
\be\label{breather}
B(t,x;\alpha,\beta,x_1,x_2) \;=\; 2\sqrt{2}\,\partial_x\!\left[\arctan\!\left(\frac{\beta}{\alpha}\,\frac{\sin(\alpha x+\delta t+x_1)}{\cosh(\beta x+\gamma t+x_2)}\right)\right],
\ee
with $\gamma=\beta(3\alpha^2-\beta^2)$ and $\delta=\alpha(\alpha^2-3\beta^2)$. At $t=0$, $x_1=x_2=0$, the profile is even, with a central peak, oscillating side lobes, and exponential envelope decay at rate $\beta$. The conserved $L^2$ mass is $\tfrac12\int B^2\,dx=4\beta$ independent of $\alpha$, so $L^2$-normalization fixes $\beta=1$ (up to the Airy rescaling $u(t,x)\mapsto\lambda u(\lambda^3 t,\lambda x)$). As $\alpha\to\infty$, one has the formal asymptotics
\be\label{breather-limit}
B(0,x;\alpha,1,0,0) \;\sim\; 2\sqrt{2}\cos(\alpha x)\,\sech(x) \;+\; O(1/\alpha),
\ee
i.e., a rapidly oscillating wave packet modulated by a fixed $\sech$ envelope. This family, and this asymptotic limit, will play a central role in what follows. Although breathers are exact nonlinear solutions of mKdV and not of the linear Airy equation, their initial data at $t=0$ is a legitimate $L^2(\R)$ function and can be used as input to the Airy--Strichartz variational problem, where their ratio $R[B(0,\cdot;\alpha,1,0,0)]$ turns out to be quantitatively meaningful.

\section{AI-based discovery pipeline}
\label{sec:method}

We now describe the AI-based numerical procedure used to prove the variational problem \eqref{sharp}. The approach has been carefully designed to match the analytic setting, both in terms of function spaces and in terms of rigorous quadrature-type approximations we have proved in \cite{MV,ACFMV,FMV,MMP}.

Most experiments were conducted on a 64-bit MacBook Pro equipped with an M4 chip, using MPS and single-precision (float32). However, when computing mixed norms of the form $L_t^qL_x^r$ with large exponents, numerical instabilities (e.g., overflow leading to NaNs) were observed in single precision. In these cases, computations were carried out in double precision (float64) on Google Colab Pro using an NVIDIA T4 GPU to ensure numerical stability. We observed no significant qualitative differences in the regimes where both precisions are stable (See Appendix \ref{app:hyper}).

\paragraph{Ansatz.} We represent the candidate initial datum as a neural network $\phi_\theta:\R^d\to\R^k$ with trainable parameters $\theta$. The architecture is fully-connected, $4$ hidden layers, $20$ neurons per layer, with the wavelet activation
\be\label{act}
\sigma(x)=e^{-(s_0x)^2}\sin(w_0 x+b_0),
\ee
with $(s_0,w_0,b_0)$ trainable and initialized at $(1/\sqrt{2},1,\mathrm{Unif}(-\pi/2,\pi/2))$. Plain $\tanh$ activations yield the same qualitative conclusions at slightly slower convergence rates. Where symmetry is expected (Schr\"odinger in $d=1$, even Airy), we symmetrize $\phi_\theta(x)\mapsto\tfrac12(\phi_\theta(x)+\phi_\theta(-x))$ before every forward pass.

\paragraph{Propagator and norm approximations.} For any of the linear propagators considered (Schr\"odinger $e^{it\Delta}$, Airy $e^{-t\partial_x^3}$, etc.), $V(t)u_0$ admits a Fourier-symbol representation
$V(t)u_0 \;=\; \F^{-1}_{\xi\to x}\!\big(\widehat{V_t}(\xi)\widehat{u_0}(\xi)\big)$, which we compute via the FFT on a uniform grid $[-R,R]_x\times[-T,T]_t$. The FFT outcome is denoted $\mathcal V[\phi_\theta]$. Mixed norms are approximated by trapezoidal quadratures $\mathcal J_X[\cdot]$; for $X=L^q_t L^r_x$, for instance,
\be\label{JLq}
\mathcal J_{L^q_tL^r_x}[v] \;=\; \Big(\Delta t\,\sum_\ell\Big(\Delta x\,\sum_j|v(t_\ell,x_j)|^r\Big)^{q/r}\Big)^{1/q},
\ee
and fractional derivatives $|D_x|^\gamma$ are applied in the frequency domain before the inverse FFT. We have proved in \cite{FMV,ACFMV,MMP} that both quantities converge rigorously to their continuous counterparts as the grid refines.

\paragraph{Loss and $L^2$ constraint.} Because the constraint $\|u_0\|_{L^2_x}=1$ is \emph{hard}---a soft penalty lets the network escape to infinity---we enforce it by pre-normalization: $\phi_\theta\leftarrow\phi_\theta/\mathcal J_{L^2_x}[\phi_\theta]$ at every forward pass. The loss is then
\be\label{loss}
\mathcal L(\theta) \;=\; -\mathcal J_X[\mathcal V[\phi_\theta]] + \mathcal L_b(\theta) + \lambda\|\theta\|_{L^2}^2,
\ee
with a mild ridge $\lambda=10^{-6}$ and a problem-specific boundary term $\mathcal L_b$; in the Airy case, for instance, $\mathcal L_b(\theta)=|\phi_\theta(0)-1|^2$ just fixes the amplitude at the origin. Optimization is done by Adam with cyclic learning rate schedule; see Appendix~\ref{app:hyper} for detailed hyperparameters.

\begin{figure}[t]
\centering
\includegraphics[width=0.8\textwidth]{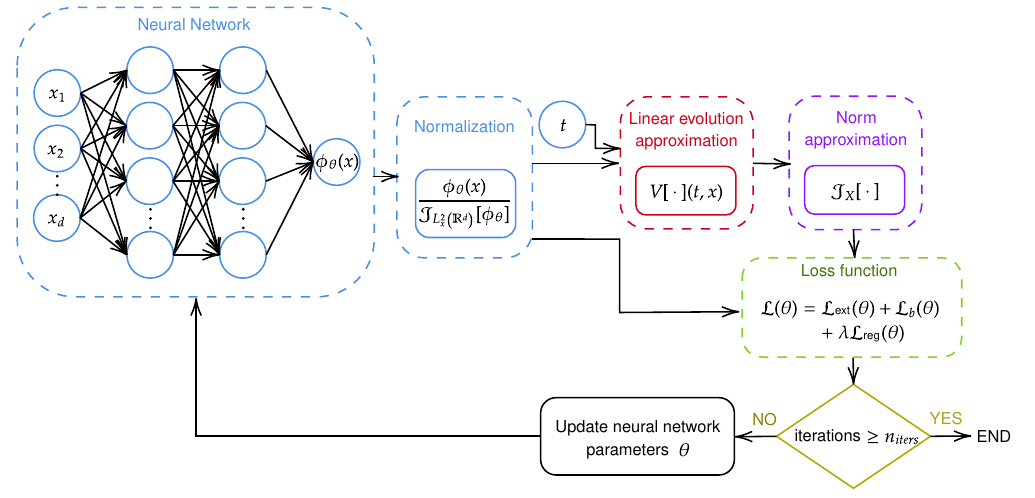}
\caption{Neural optimization pipeline. The network $\phi_\theta$ is normalized in $L^2_x$ at every forward pass, propagated by the FFT-based linear group $\mathcal V[\cdot]$, and evaluated in the target mixed norm $\mathcal J_X$ to produce the loss~\eqref{loss}.}
\label{fig:arch}
\end{figure}

\paragraph{What convergence means.} We declare that the neural optimization has found a candidate extremizer if, as training proceeds, (a) $\phi_\theta$ stabilizes at a bounded-$L^2$ profile and (b) $\mathcal J_X[\mathcal V[\phi_\theta]]$ converges to a finite value. Failure of (a), manifested by progressive scale expansion or oscillatory drift of $\phi_\theta$, is interpreted as numerical evidence of the non-existence of extremizers. This is the regime we will encounter in the Airy--Strichartz case below.

\section{Validation on the linear Schr\"odinger model}
\label{sec:schrod}

\begin{figure}[t]
\centering
\begin{subfigure}[t]{0.45\linewidth}\centering
\includegraphics[width=\linewidth]{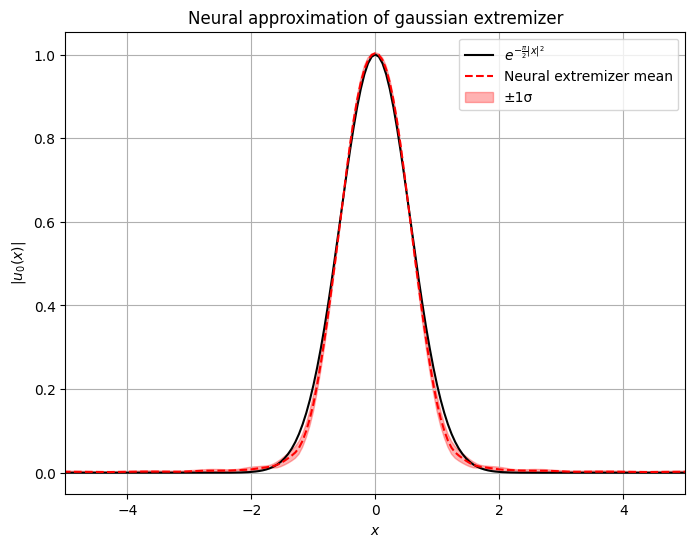}
\caption{Neural extremizer $\phi_{\theta^\star}$ versus analytical Gaussian.}
\end{subfigure}\hfill
\begin{subfigure}[t]{0.45\linewidth}\centering
\includegraphics[width=\linewidth]{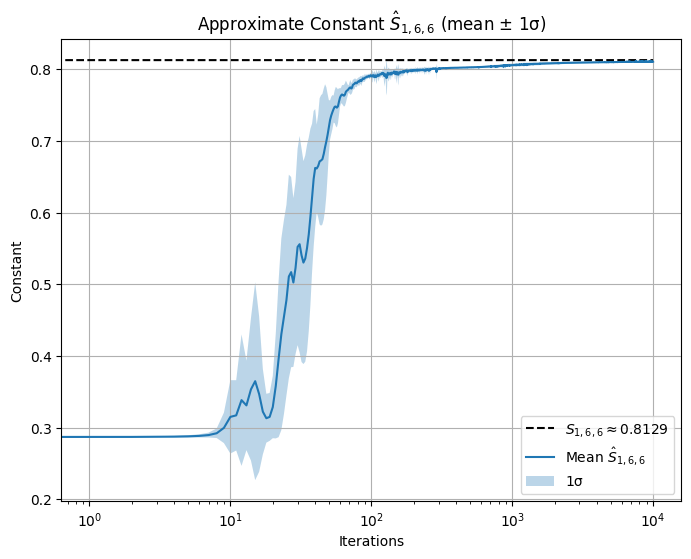}
\caption{Evolution of $\hat S_{1,6,6}$ during training.}
\end{subfigure}
\caption{Schr\"odinger $d=1$, $p=6$: the neural network recovers Foschi's Gaussian extremizer with no analytical input.}
\label{fig:schrod}
\end{figure}

 The Schr\"odinger case provides a clean double test: for $(d,p)\in\{(1,6),(2,4)\}$ the extremizer is a Gaussian \cite{Kunze03,Foschi07,HundertmarkZharnitsky06}, while for other admissible pairs with $q\ne r$ the extremizer exists \cite{Shao09} but its shape is still unknown. The pipeline should find Gaussians in the first case and identify the shape of the extremizer in the second case. Figure~\ref{fig:schrod} shows the $d=1$ results: the neural-approximated profile is graphically indistinguishable from the analytical Gaussian $u_0^\star\propto e^{-bx^2}$ of Foschi, and the discovered constant $\hat S_{1,6,6}$ matches Foschi's value to relative error $\lesssim 10^{-3}$. See Appendix \ref{sec:3} for further validation on Schr\"odinger, including a table reporting $\widehat S_{1,6,6}$ and $\widehat S_{2,4,4}$ vs their analytical values, with error estimates, for $d=1,2$.

\subsection{Results in the non-endpoint Schr\"odinger admissible pairs}
\begin{figure}[t]
\centering
\begin{subfigure}[t]{0.45\linewidth}\centering
\includegraphics[width=\linewidth]{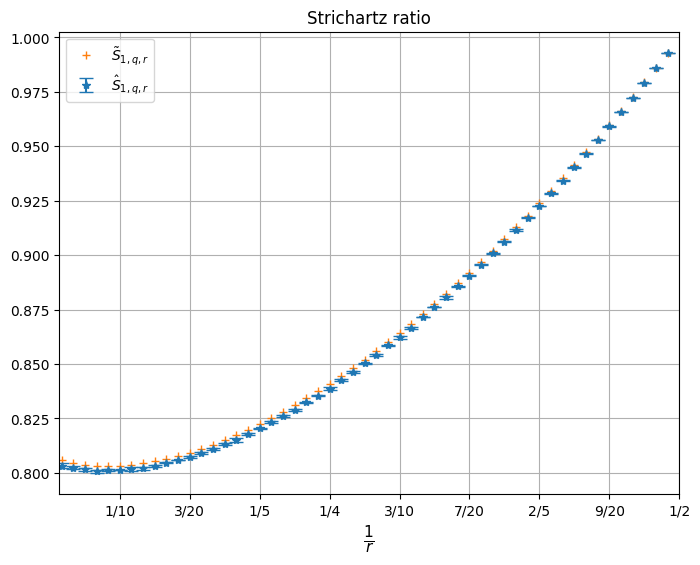}
\caption{Neural extremizer $\phi_{\theta^\star}$ ratio versus analytical gaussian ratio.}
\end{subfigure}\hfill
\begin{subfigure}[t]{0.45\linewidth}\centering
\includegraphics[width=\linewidth]{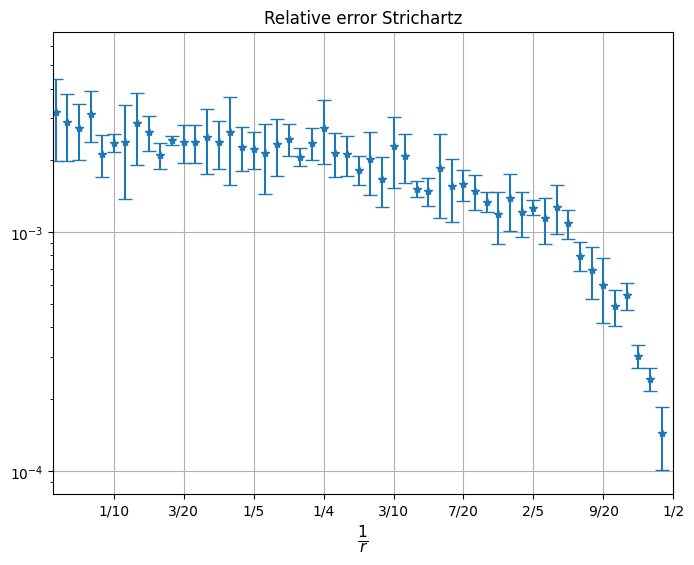}
\caption{Relative error of  $\hat{S}_{1,q,r}$, with respect to the gaussian ratio.}
\end{subfigure}
\caption{Schr\"odinger $d=1$: the neural optimization finds Gaussians as extremizers for any admissible pair.}
\label{fig:schrod2}
\end{figure}

 We run the full neural optimization pipeline for $59$ Schr\"odinger-admissible pairs $(q,r)$. Figure~\ref{fig:schrod2} shows the approximate value of the Strichartz ratio via NNs, $\hat{S}_{1,q,r}$, and the relative error with respect to the Strichartz ratio of a Gaussian, $\widetilde S_{1,q,r}$. Error bars are obtained from $5$ independent runs of the algorithm for each pair and report the standard deviation. As can be seen in the figure, the neural-approximated profile corresponds to a Gaussian as well, with a relative error $\lesssim 10^{-3}$ for all the admissible pairs. Each training run takes less than 3 minutes of wall-clock time.
 
\section{Airy--Strichartz: main findings}
\label{sec:airy}

We now turn to the critical Airy--Strichartz inequality \eqref{AS} with $\gamma=1/q$, which is the main focus of this paper.

\subsection{Baseline: soliton and breather ratios}

For any $w\in L^2(\R)$ we write
\[
R[w] := \frac{\|\,|D_x|^\gamma e^{-t\partial_x^3}w\|_{L^q_tL^r_x}}{\|w\|_{L^2_x}}.
\]
By the Airy scaling $u(t,x)\mapsto\lambda u(\lambda^3 t,\lambda x)$, $R[Q_{p,c}]$ is independent of $c$ and $R[B(0,\cdot;\alpha,\beta,\cdot,\cdot)]$ depends only on $\alpha/\beta$. The effective parameter spaces are therefore the discrete $p$ (for solitons) and the continuous $\alpha$ (for breathers with $\beta=1$).

\begin{figure}[t]
\centering
\begin{subfigure}[t]{0.45\linewidth}\centering
\includegraphics[width=\linewidth]{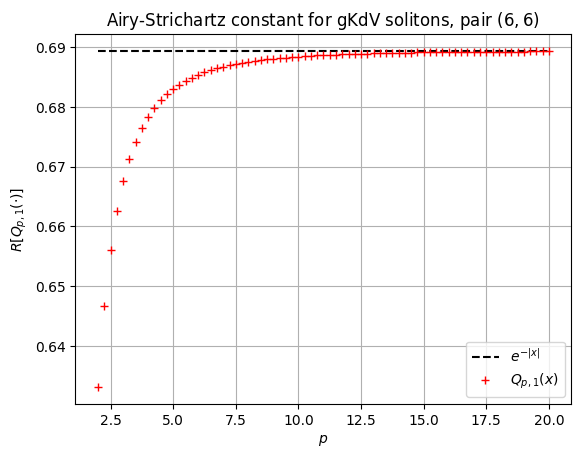}
\caption{Soliton ratios $R[Q_{p,1}]$ (red) and limit $R[e^{-|x|}]$ (black) at $(q,r)=(6,6)$.}
\label{fig:ratio-sol}
\end{subfigure}\hfill
\begin{subfigure}[t]{0.45\linewidth}\centering
\includegraphics[width=\linewidth]{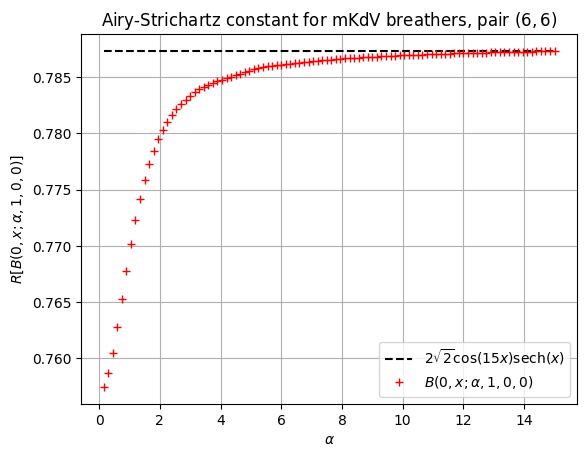}
\caption{Breather ratios $R[B(0,\cdot;\alpha,1,0,0)]$ (red) and limit $R[2\sqrt{2}\cos(\alpha x)\sech(x)]$ (black).}
\label{fig:ratio-breather}
\end{subfigure}
\caption{Airy--Strichartz ratios for solitons and breathers at $(q,r)=(6,6)$, in both cases below $\widetilde A_{6,6}$.}
\end{figure}

Figures~\ref{fig:ratio-sol}--\ref{fig:ratio-breather} summarize two facts. First, $R[Q_{p,1}]$ is monotonically increasing in $p$ and is bounded by the ratio of the formal scaling limit $R[e^{-|x|}]$: solitons are not maximizing, and they are not even close. Second, $R[B(0,\cdot;\alpha,1,0,0)]$ is monotonically increasing in $\alpha$ and approaches, but does not cross, the limit value $R[2\sqrt{2}\cos(\alpha x)\sech(x)]$, which in turn lies very close to $\widetilde A_{6,6}$. Breathers therefore exceed every soliton ratio and approach the Frank--Sabin bound from below---but only asymptotically, as $\alpha\to\infty$. Already this baseline suggests that if anything saturates the Frank--Sabin bound for \eqref{AS}, it is not a single $L^2$ profile but a sequence.

\subsection{Hermite ansatz: an independent check}

As a non-neural cross-check, we fit a $21$-term Hermite expansion $u_0(x)=\sum_{n=0}^{20}b_n f_n(x)$ and maximize $R[u_0]$ over the coefficients. The resulting profile is visually and quantitatively very close to an mKdV breather $B(t_0,\cdot;\alpha,\beta,0,x_0)$ with $\alpha\approx 1.24$, $\beta\approx 0.63$, giving a ratio $\hat A_{6,6}^{\text{Hermite}}=0.7828<\widetilde A_{6,6}\approx 0.7886$. The Hermite basis is orthonormal in $L^2(\R)$ and has no built-in knowledge of breathers, so its independent convergence to a breather is non-trivial. The coefficients, the truncation analysis (the answer stabilises only for $N\gtrsim 18$, with large cancellations among modes $7$--$13$) and the figure are reported in Appendix~\ref{app:hermite}.

\subsection{Numerical results and non-existence evidence}

\begin{figure}[t]
\centering
\begin{subfigure}[t]{0.45\linewidth}\centering
\includegraphics[width=\linewidth]{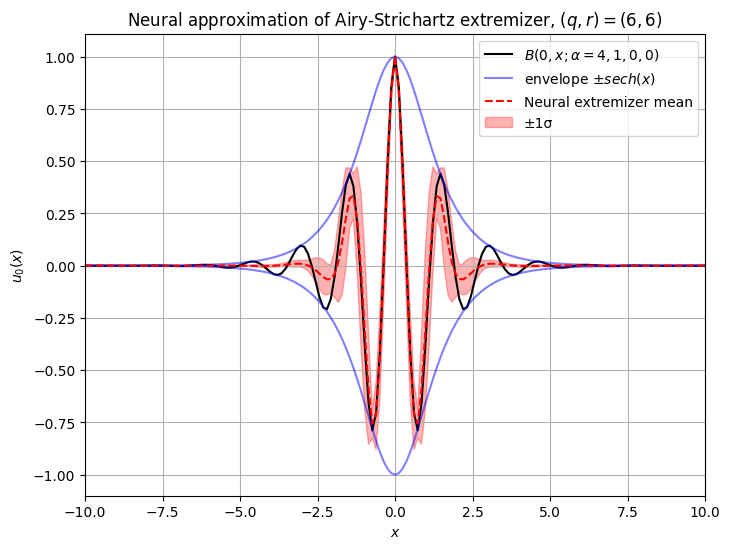}
\caption{Approximate extremizer $\phi_{\theta^\star}(x)$ for the pair $(6,6)$.}
\end{subfigure}\hfill
\begin{subfigure}[t]{0.45\linewidth}\centering
\includegraphics[width=\linewidth]{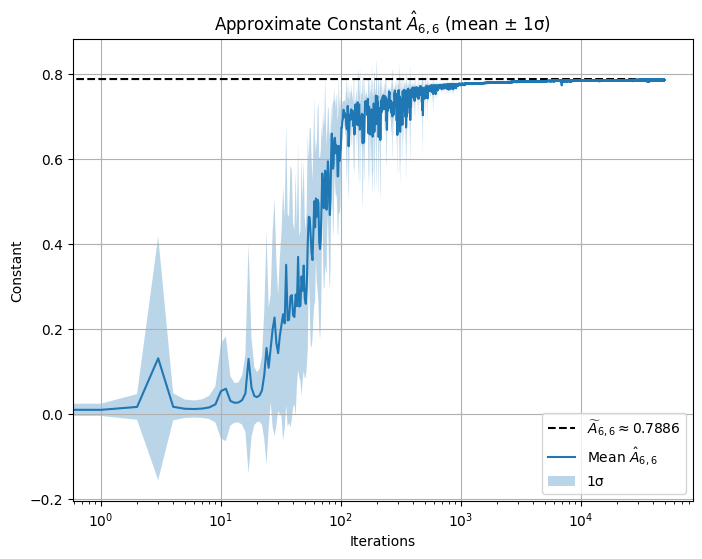}
\caption{Convergence of $\hat A_{6,6}$ towards $\widetilde A_{6,6}$.}
\end{subfigure}
\caption{AI-based discovery for the Airy--Strichartz inequality, pair $(6,6)$.}
\label{fig:pinn-airy}
\end{figure}

We run the full neural optimization pipeline of Section~\ref{sec:method} on the four admissible pairs $(q,r)\in\{(5,10),(6,6),(8,4),(12,3)\}$, all at $\gamma=1/q$. Figure~\ref{fig:pinn-airy} plots the discovered profile for the pair $(6,6)$ in comparison with an mKdV breather with $\alpha=4$ and the convergence of the discovered ratio $\hat A_{6,6}$ with training iterations. In this setting, we obtain a relative error of $2.576 \times 10^{-3}$ (against the Frank--Sabin lower bound $\widetilde A_{6,6}\approx 0.7886$) with the estimated constant $\hat{A}_{6,6} \approx 0.7865$. More details on the results for the other pairs, including the relative errors, can be found in Appendix \ref{app:AS-numerics}. Each training run takes approximately 45 minutes of wall-clock time.

\begin{table}[t]
\centering
\caption{Evolution of $\hat A_{6,6}$ and the error (mean $\pm$ std) across iterations, compared with the Frank--Sabin lower bound $\widetilde A_{6,6}$.}
\label{tab:A66}
\begin{tabular}{lcc}
\toprule
Iterations & $\hat A_{6,6}$ & ${\bf error}_{6,6}$ \\
\midrule
100   & $0.6168 \pm 0.1311$ & $(2.179 \pm 1.662)\times 10^{-1}$ \\
1000    & $0.7773 \pm 0.0022$ & $(1.436 \pm 0.283)\times 10^{-2}$ \\
10000   & $0.7852 \pm 0.0010$ & $(4.345 \pm 1.317)\times 10^{-3}$ \\
50000   & $0.7865 \pm 0.0011$ & $(2.576 \pm 1.363)\times 10^{-3}$ \\
\bottomrule
\end{tabular}
\end{table}

\begin{table}[t]
\centering
\caption{Consolidated NN results across the four critical Airy--Strichartz pairs at $\gamma=1/q$, after $5\times 10^4$ iterations. Mean over $5$ independent runs; relative error is computed against the Frank--Sabin lower bound $\widetilde A_{q,r}$. In every pair $\hat A_{q,r}<\widetilde A_{q,r}$ by a shrinking margin, consistent with non-attainment.}
\label{tab:Aqr-summary}
\begin{tabular}{lccc}
\toprule
$(q,r)$ & $\widetilde A_{q,r}$ & $\hat A_{q,r}$ & Rel.\ error \\
\midrule
$(5,10)$  & $0.7926$ & $0.7915$ & $1.642\times 10^{-3}$ \\
$(6,6)$   & $0.7886$ & $0.7865$ & $2.576\times 10^{-3}$ \\
$(8,4)$   & $0.8112$ & $0.8084$ & $3.385\times 10^{-3}$ \\
$(12,3)$  & $0.8556$ & $0.8541$ & $1.680\times 10^{-3}$ \\
\bottomrule
\end{tabular}
\end{table}

Three observations are in order. \emph{(i) No $L^2$ stabilization.} Unlike the Schr\"odinger case, where $\phi_\theta$ stabilizes at a fixed Gaussian, in the Airy case the neural profile continues to gain internal oscillations as training proceeds, i.e., its effective frequency parameter $\alpha_{\max}$ grows with the number of iterations. No fixed $L^2$ profile is selected by the training dynamics; what stabilizes is the \emph{shape of the partially-truncated breather} at the current $\alpha_{\max}$. \emph{(ii) Quantitative breather fit.} Figure~\ref{fig:pinn-airy} shows that, at the final training snapshot, the neural profile fits an mKdV breather $B(0,\cdot;\alpha,1,0,0)$ quantitatively, with the decay envelope $|u_0(x)|\lesssim \sech(x)$, essentially matching the breather with $\beta=1$. This rules out polynomial or Gaussian-tailed alternatives, which would produce a very different log-slope. \emph{(iii) Approach to the Frank--Sabin bound from below.} The ratio $\hat A_{q,r}$ approaches $\widetilde A_{q,r}$ monotonically but always stays strictly below it, with the remaining gap shrinking as training progresses. This is the expected signature of $A_{q,r}=\widetilde A_{q,r}$: the supremum is realized only as the weak limit of a non-compact maximizing sequence, and a neural network with finite grid and finite network size can only approximate the $\alpha\to\infty$ member of the sequence up to some effective $\alpha_{\max}$. The Frank--Sabin theorem then explicitly predicts that the bound will not be crossed by any compact approximation.

\subsection{Airy-Strichartz for a non critical case}
\begin{figure}[!ht]
\centering
\includegraphics[width=0.45\linewidth]{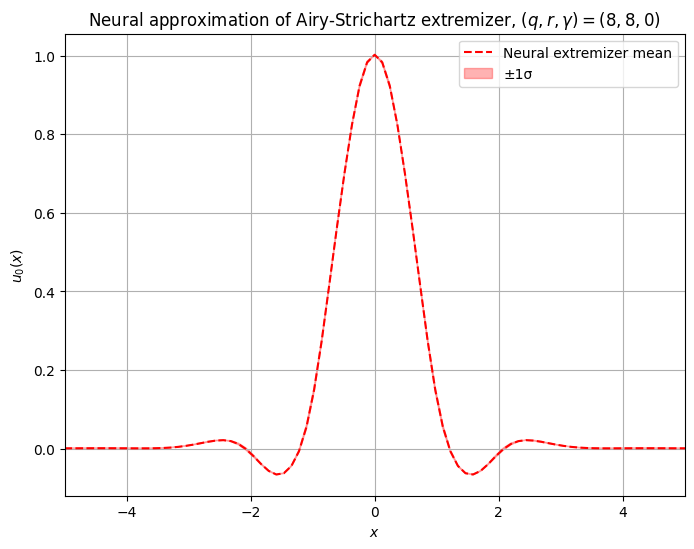}
\caption{Approximate extremizer $\phi_{\theta^\star}(x)$ for the pair $(8,8)$.}
\label{fig:A88}
\end{figure}

Besides the critical case $\gamma=\frac1q$, there is a non-critical case where extremizers exist. Indeed, for the pair $(q,r)=(8,8)$ we have $\gamma = 0$, and extremizers for the Airy-Strichartz inequality are well known to exist \cite{HundertmarkShao10}, but the extremizer profile is still unknown. We perform the AI-based algorithm to find the approximated extremizer profile. Figure~\ref{fig:A88} shows the neural profile, for 10 independent seeds. A fitted $L^2$ profile for this neural extremizer, as well as information regarding $\hat A_{8,8}$ can be found in Appendix~\ref{app:AS-noncrit}. We have found a possible candidate for extremizer in a class of modified gaussians.

\subsection{Main conjecture}

Combining the soliton/breather baseline, the independent Hermite check, and the NN results, we arrive at the following.

\begin{conjecture}\label{conj:main}
For every Schr\"odinger-admissible pair $(q,r)$ with $4<q<\infty$ and $\gamma=1/q$, the sharp Airy--Strichartz constant defined in \eqref{Aqr} satisfies $A_{q,r}=\widetilde A_{q,r}$, where $\widetilde A_{q,r}$ is the Frank--Sabin bound \eqref{FS-bound}. Moreover, the supremum is not attained by any $u_0\in L^2(\R)$, but is realized, modulo the symmetries of the Airy equation, as the limit
\be\label{limit-breather}
A_{q,r} \;=\; \lim_{\alpha\to\infty}\, R[B(0,\cdot;\alpha,1,0,0)],
\ee
with $B$ the mKdV breather \eqref{breather}.
\end{conjecture}

Conjecture~\ref{conj:main} fits the Frank--Sabin dichotomy exactly: strict inequality $A_{q,r}>\widetilde A_{q,r}$ would produce an $L^2$ extremizer by precompactness, whereas equality $A_{q,r}=\widetilde A_{q,r}$ is precisely the regime in which the supremum is realized only as a non-compact limit. Our numerical evidence is consistent with the second scenario. A proof would plausibly rely on a direct evaluation of $R[B(0,\cdot;\alpha,1,0,0)]$ in the limit $\alpha\to\infty$ together with the sharp Schr\"odinger constant $S_{1,q,r}$ and the identities of \cite{HundertmarkZharnitsky06,BBC H09}; we outline the strategy in Appendix~\ref{app:outlook}.

\section{Discussion}
\label{sec:discussion}

\subsection{Conclusions}

The Schr\"odinger experiments confirm what was already known from \cite{Foschi07,HundertmarkZharnitsky06} at the sharp exponents and provide consistent computational support for the conjecture that Gaussians remain the extremizers across the entire admissible range. The novelty lies in the Airy case: two independent ans\"atze (Neural optimization and Hermite) converge to the same family of mKdV breathers, the discovered ratio approaches the Frank--Sabin lower bound monotonically from below, and the residual gap shrinks as a clean power law $\sim\alpha^{-0.9}$. Conjecture~\ref{conj:main} crystallises this picture.

The reason why mKdV breathers appear here is, in retrospect, natural. The Airy group is the linear part of focusing mKdV, and breathers are the coherent structures that most economically trade off spatial concentration and temporal oscillation. The balance between the two parameters $(\alpha,\beta)$ is exactly the balance a scale-invariant space-time inequality rewards. As $\alpha\to\infty$ the breather becomes an infinitely oscillatory $\sech$-enveloped wave packet, whose Strichartz ratio approaches but never attains the Frank--Sabin bound; this, in our reading, is why no explicit $L^2$ extremizer for the critical Airy--Strichartz inequality has ever been proposed.

The pipeline applies verbatim to the Wave, Klein--Gordon and Zakharov--Kuznetsov groups, where existence of extremizers and the sharpness of Frank--Sabin-type criteria remain open. We view the breather-as-asymptote picture as a template that may be revisited analytically in those settings. Operator-learning architectures (DeepONet, FNO) could complement the present approach, but they do not directly produce extremizers; we leave their use as pretraining mechanisms for future work.

\subsection{Limitations and reproducibility}

The empirical core of our non-existence evidence is the ``growing-$\alpha$'' signature, observed up to numerically tractable $\alpha$; we have run additional experiments at higher $\alpha$ with the same monotone behaviour, but the cost of resolving the high-frequency oscillations grows fast. The non-existence claim itself is delicate: numerical evidence cannot strictly distinguish non-attainment from very slow attainment. What we provide is a consistent picture across two independent ans\"atze together with the quantitative power-law fit $\widetilde A_{6,6} - R[B(0,\cdot;\alpha,1,0,0)] \sim C\alpha^{-0.9}$ documented in Appendix~\ref{app:estimates}; we view this as a strong invitation to a rigorous proof along the analytical route outlined in Appendix~\ref{app:outlook}, not as a substitute for one. Finally, we have restricted attention to the Schr\"odinger and Airy groups; the wave, Klein--Gordon and Zakharov--Kuznetsov cases are left for future work.

The full code (NN training scripts, FFT-based propagators, Hermite-ansatz solver, post-processing notebooks for tables and figures) together with trained checkpoints will be released anonymously at submission time. Exact hyperparameters are reported in Appendix~\ref{app:hyper}.

\subsubsection*{Acknowledgments}
We deeply thank Emanuel Carneiro for explaining the problem to us and for very useful comments. Part of this work was done while the third author was visiting ICTP in December 2025 as part of the Ramanujan Prize 2025 ceremony; he acknowledges the support and warm hospitality received. R.F.\ was partially funded by Chilean research grants FONDECYT 3230256, MathAmSud WAFFLE and ANID Exploraci\'on 13220060. C.M.\ was partially funded by Chilean research grants FONDECYT 1231250, MathAmSud WAFFLE, ANID Exploraci\'on 13220060 and Basal CMM FB210005. N.V.\ was partially funded by Chilean research grants ANID 2022 Exploration 13220060, FONDECYT 1231250, a Latin America PhD Google Fellowship, ANID-Subdirecci\'on de Capital Humano/Doctorado Nacional/2023-21231021 and Basal CMM FB210005. All authors were funded by the INRIA PANDA innitiative, whose support is greatly acknowledged.

\bibliographystyle{plainnat}

\appendix

%
%

\section{The linear Schr\"odinger case}\label{sec:3}

Recall the Schr\"odinger--Strichartz estimate \eqref{Str-Sch}, and define the sharp constant
\begin{equation}\label{Sdqp}
S_{d,q,r}\;:=\;\sup\bigl\{\|e^{it\Delta}u_0\|_{L_t^{q}L_x^{r}(\R\times\R^d)}:\ \|u_0\|_{L^2(\R^d)}=1\bigr\}.
\end{equation}
Foschi's Theorem and the sharp constants for $(d,q,r)\in\{(1,6,6),(2,4,4),(1,8,4)\}$ etc., proved in \cite{Foschi07,HundertmarkZharnitsky06,Carneiro09,BBC H09}, are summarised in Table~\ref{tab:Strichartz-ctes}. We will also use, in the validation tables below, the following normalization observation.
\begin{remark}\label{rem:unit}
Under \eqref{Str-Sch} and the assumption $|u_0(0)|=1$ together with $L^2$ unit mass and even parity, Foschi's Gaussian extremizer $\widetilde f(x)=\exp(A|x|^2+b\cdot x+C)$ reduces, modulo phase, to $|u_0(x)|=e^{-\frac{\pi}{2}|x|^2}$.
\end{remark}

\begin{table}[!ht]
\centering
\begin{tabular}{|c|c|c|c|}
\hline
$(d,q,r)$ & Type & Optimal constant \eqref{Sdqp} & Reference \\
\hline
$(1,6,6)$ & Strichartz & $12^{-1/12}$ & \cite{Foschi07} \\\hline
$(2,4,4)$& Strichartz  & $2^{-1/2}$ & \cite{Foschi07} \\\hline
$(1,8,4)$& Strichartz  & $2^{-1/4}$ & \cite{HundertmarkZharnitsky06} \\\hline
$(1,12,6)$ & Sobolev-Strichartz  & $(6\pi)^{-1/12}$        & \cite{Carneiro09} \\\hline
$(1,16,4)$ & Sobolev-Strichartz  & $(8\pi)^{-1/16}$        & \cite{Carneiro09} \\\hline
$(2,6,6)$ & Sobolev-Strichartz   & $(12\pi)^{-1/6}$        & \cite{Carneiro09} \\\hline
$(2,8,4)$ & Sobolev-Strichartz   & $(16\pi)^{-1/8}$        & \cite{Carneiro09} \\\hline
$(4,4,4)$ & Sobolev-Strichartz   & $(32\pi)^{-1/4}$        &  \cite{Carneiro09}\\\hline
\end{tabular}
\caption{Sharp Strichartz estimates in various dimensions $d$, with the reference where the optimal constant is established. Equality holds iff $f$ is in the Gaussian family.}
\label{tab:Strichartz-ctes}
\end{table}

\subsection{The case of soliton solutions} 
\begin{figure}[!ht]
\centering
\includegraphics[width=0.5\textwidth]{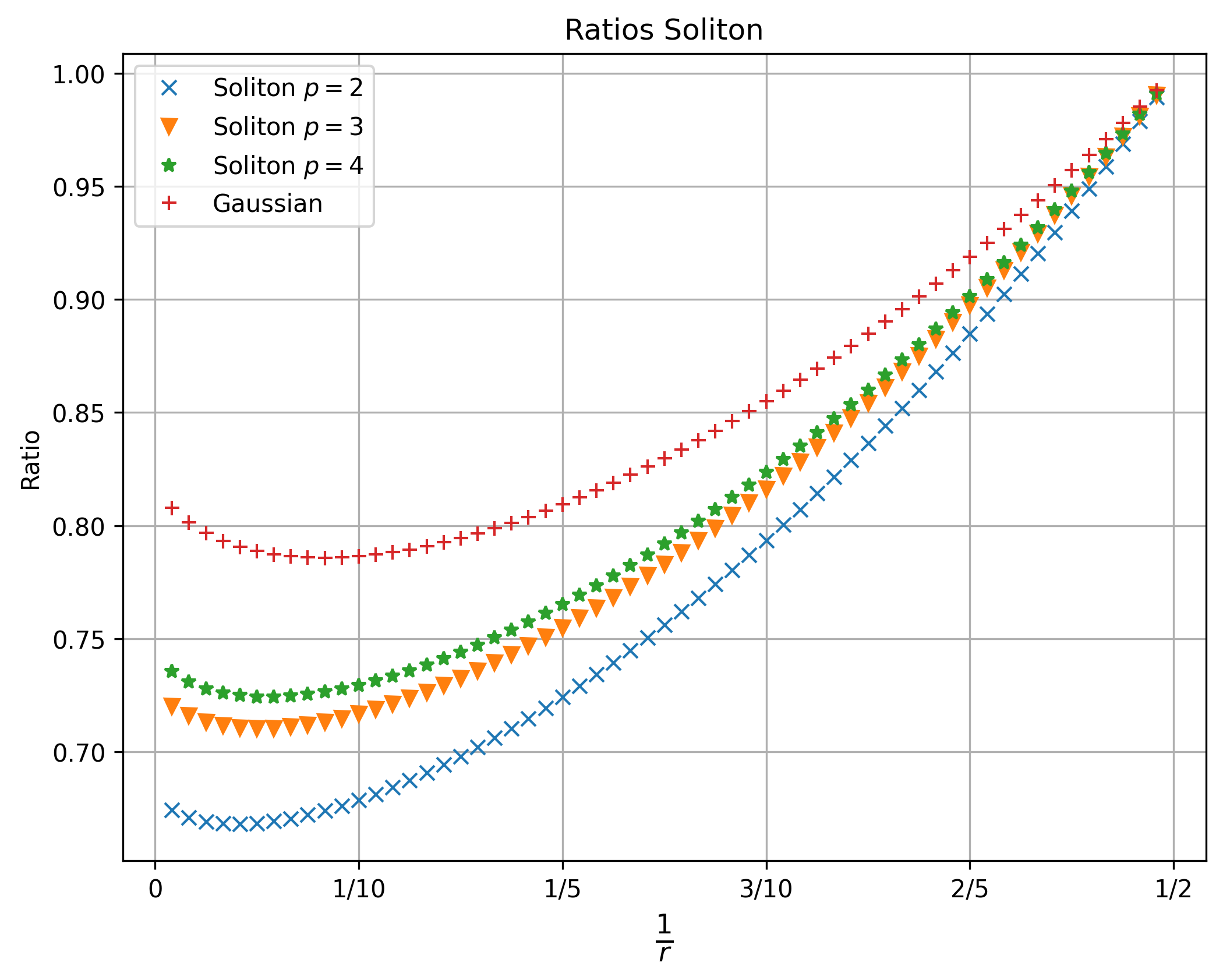}
\caption{Ratio comparison between three solitons $Q_p$, with $p \in \{2,3,4\}$, and the gaussian $e^{-\frac{\pi}{2}|x|^2}$. The comparison is made in function of $\frac 1r \in (0,\frac12)$, taking 59 different values. In the case $r=2$ we have the classical $L^2$ conservation and every $L^2$ solution has the same ratio.}
\label{fig:ratios-soliton}
\end{figure}

An important case of a possible candidate for an extremal point in the previous analysis is the soliton itself. Even if we know that in some cases Gaussians extremize Strichartz estimates, it is also important to regard in what sense solitons are not extremizers, or how far they are from being extremizers. To fix ideas, we consider the 1D case only. Consider the nonlinear Schrodinger model
\be\label{NLS}
i\partial_t u + \partial_x^2 u + |u|^{p-1} u=0, \quad p>1,
\ee
with $u=u(t,x)\in\mathbb C$ and $t,x\in\mathbb R$. The solitary waves for \eqref{NLS} are given by
\be\label{NLS_SW}
S[p,c,v,x_0,\gamma](t,x):=Q_{p,c}(x-vt-x_0) e^{ixv/2}e^{ict}e^{itv^2/4 +i\gamma},
\ee
where $Q_{p,c}(s)=c^{\frac1{p-1}}Q_p(\sqrt{c}s)$, and 
\be\label{Qk}
Q_p(s):=\left(\frac{p+1}{2\cosh^2\left(\frac{p-1}2 s\right)}\right)^{\frac1{p-1}}.
\ee
Let $(d=1,q,r)$ be an admissible triple. Notice that \eqref{NLS} preserves the $L^2$ norm, then we can consider just $t=0$. By natural shifts and scaling symmetries, we can assume $x_0=\gamma=0$ and $c=1$. In this direction, we first quickly realize from the Galilean symmetry in integration on $\mathbb R_t\times \mathbb R_x$ (in the case of the Strichartz estimate) that the ratio
\[
\frac{\|e^{t\partial_x^2} S[p,1,v,0,0](0,\cdot)\|_{L^q_tL^r_x}}{\|S[p,1,v,0,0](0,\cdot)\|_{L^2}} = \frac{\|e^{t\partial_x^2} (Q_{p,1}(x)e^{ixv/2})\|_{L^q_tL^r_x}}{\|Q_{p,1}\|_{L^2}} =: R[p,q].
\]
is also independent of $v$. We now present numerical approximations of the ratio $R[p,q]$ for different values of $p$ and $q$. The results are summarised in Figure~\ref{fig:ratios-soliton}, which shows the ratios for $59$ values of $\frac1r$ on a uniform grid over $(0,\tfrac12)$ with $\Delta r=\tfrac1{120}$, for initial data given by three different solitons and the gaussian $e^{-\frac{\pi}{2}|x|^2}$. The gaussian attains the largest ratio; the ratio for solitons decreases as $p$ decreases. All four initial conditions approach $1$ as $\frac1r\to 0$, as a consequence of mass conservation.

An important remark for gaussians: they enjoy shift and scaling symmetries, so the ratio of any gaussian equals the ratio of the $L^2$-normalised gaussian $e^{-\frac{\pi}{2}|x|^2}$. A direct computation gives
\be\label{tildeS}
R\!\left[e^{-\frac{\pi}{2}|\cdot|^2}\right] = \left(\frac 2r\right)^{\frac d{2r}} \left(\frac12\right)^{\frac1q} =: \widetilde S_{d,q,r}.
\ee


\subsection{Numerics}
 We follow the neural network architecture, $L^2$-normalization, and loss defined in Section~\ref{sec:method}; we record below the full validation tables.
We shall work in dimensions $d=1,2$. In both dimensions, a neural network with 4 hidden layers and 20 neurons per hidden layer will be considered. The domain boundaries, and the number of points in the grid for the definition of $\mathcal J_{L_x^2(\R^d)}$ and $\mathcal J_{L_t^q L_x^r(\R^{d+1})}$ will be declared in each sub-case.  Additionally, the approximate constant is given as
\be\label{Sdpq-appr}
\hat{S}_{d,q,r} := -\left(- \left(\mathcal J_{L_t^qL_x^r(\R^{d+1})}[\mathcal V[\phi_{\theta^*}]]\right)^q\right)^{1/q},
\ee
where we recall $\theta^*$ to be the optimal parameters obtained in the minimization procedure. The estimate \eqref{Sdpq-appr}  is later compared with the exact value of $S_{d,q,r}$, given in Table \ref{tab:Strichartz-ctes}, by using the relative error
\[
{\bf error}_{\text{rel},d,q,r} = \frac{\big|S_{d,q,r} - \hat{S}_{d,q,r}\big|}{S_{d,q,r}}.
\]
For pairs with unknown value of $S_{d,q,r}$, we will use the conjectured value $\widetilde S_{d,q,r}$ instead, defined in \eqref{tildeS}. Finally, the absolute value $|\phi_{\theta^*}(x)|$ is compared with the gaussian function $e^{-\frac{\pi}{2}|x|^2}$, following Remark \ref{rem:unit}. 

\subsubsection{Dimension 1. Even approximations}
\begin{figure}[!ht]
\centering
\begin{subfigure}[t]{0.45\linewidth}\centering
\includegraphics[width=\linewidth]{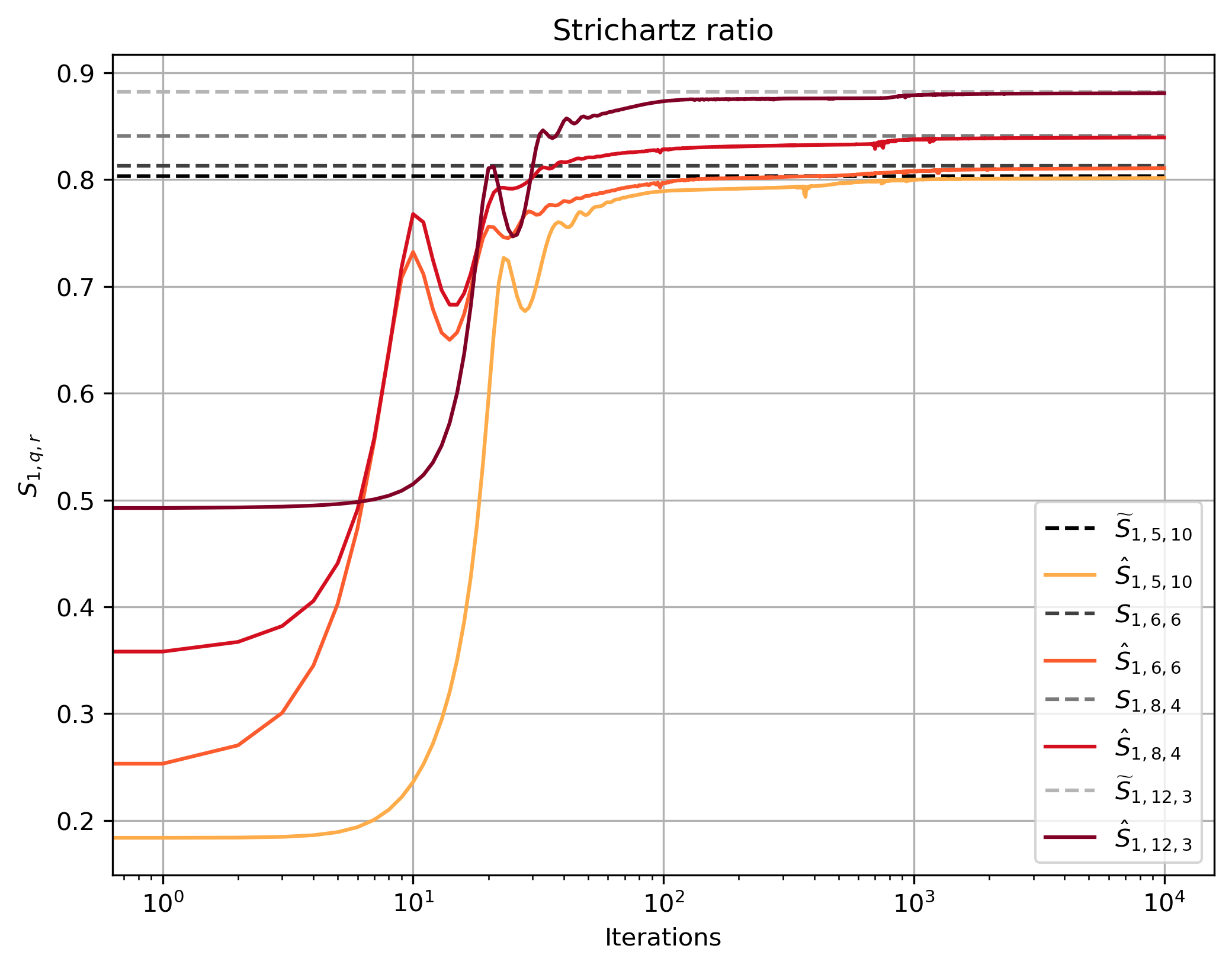}
\caption{Approximate constant $\hat{S}_{1,q,r}$.}
\end{subfigure}\hfill
\begin{subfigure}[t]{0.45\linewidth}\centering
\includegraphics[width=\linewidth]{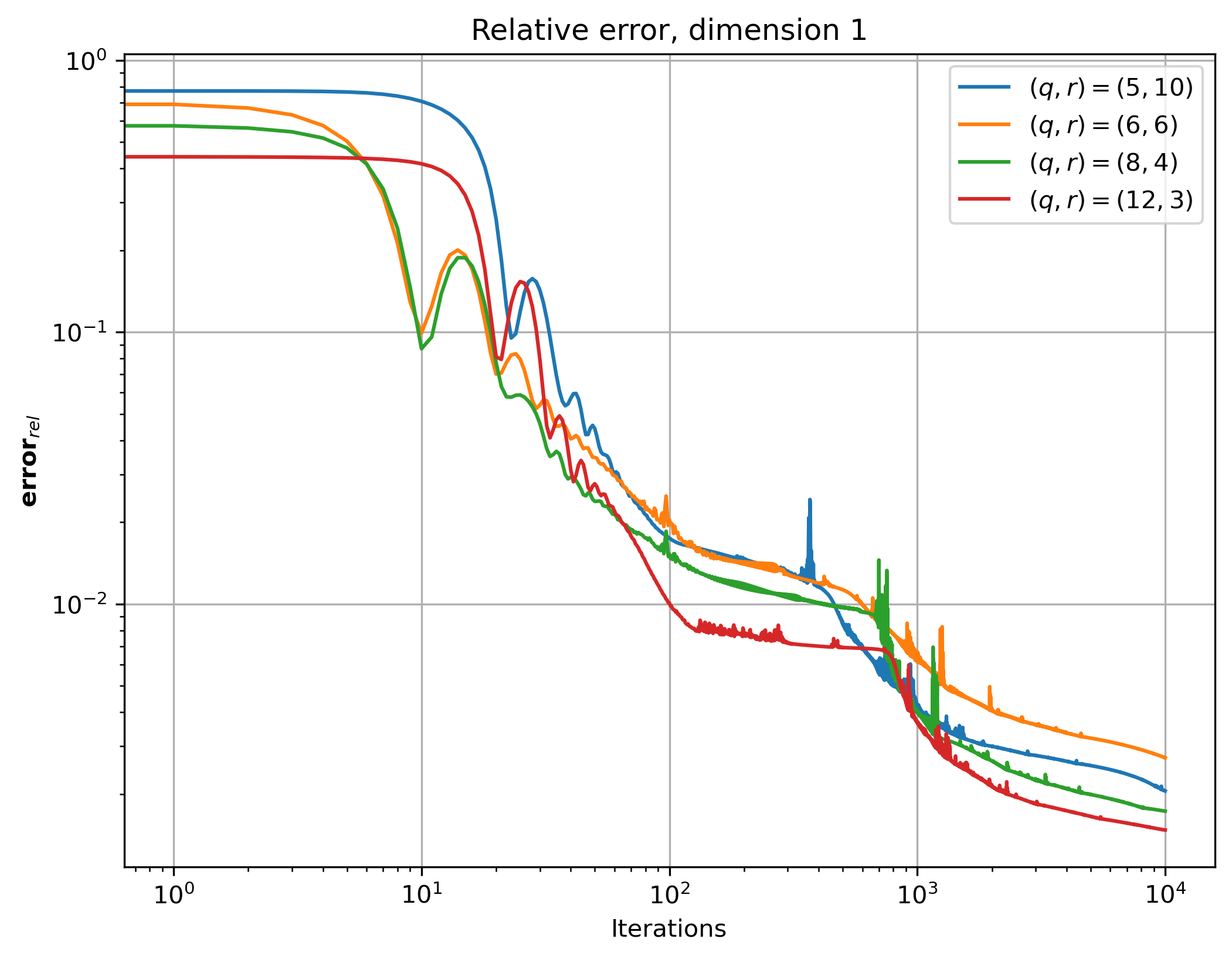}
\caption{Evolution of the relative error ${\bf error}_{\text{rel},1,q,r}$.}
\end{subfigure}
\caption{Schr\"odinger $d=1$: NN results for four different admissible pairs $(q,r)$.}
\label{fig:S1pq}
\end{figure}

In this setting, the functionals $\mathcal J_{L_x^2(\R^d)}$ and $\mathcal J_{L_t^qL_x^r(\R^{d+1})}$ are computed with $R=30$, $T=1$ and $M=N=1024$. We perform experiments in two regimes: pairs $(d,q,r)\in\{(1,6,6),(1,8,4)\}$ where the optimal constant $S_{d,q,r}$ is known, and pairs $(d,q,r)\in\{(1,5,10),(1,12,3)\}$ where the constant is conjectural.
Table~\ref{tab:S1pq} shows the evolution of the approximation $\hat S_{1,q,r}$ over iterations, averaged over $5$ independent realizations of the NN algorithm for each configuration $(q,r)$. The evolution of $\mathbf{error}_{\text{rel}}$ and $\hat S_{1,q,r}$ for a single realization of each pair is displayed in Figure~\ref{fig:S1pq}.
\begin{table}[!ht]
\def\arraystretch{1.4}
\begin{tabular}{|c||c|c|c|c|c|}
\hline
Iterations                        & 100 & 500 & 1000 & 5000 & 10000 \\ \hline \hline
$\hat S_{1,5,10}$                 & 0.7665 & 0.7922 & 0.7961 & 0.8008 & 0.8014 \\ \hline
${\bf error}_{\text{rel},1,5,10}$ & 4.573 $\times 10^{-2}$ & 1.377 $\times 10^{-2}$ & 8.904 $\times 10^{-3}$ & 3.050 $\times 10^{-3}$ & 2.319 $\times 10^{-3}$ \\ \hline \hline
$\hat S_{1,6,6}$                  & 0.7717 & 0.8008 & 0.8066 & 0.8113 & 0.8116 \\ \hline
${\bf error}_{\text{rel},1,6,6}$  & 5.070 $\times 10^{-2}$ & 1.500 $\times 10^{-2}$ & 7.777 $\times 10^{-3}$ & 1.980 $\times 10^{-3}$ & 1.626 $\times 10^{-3}$ \\ \hline \hline
$\hat S_{1,8,4}$                  & 0.8145 & 0.8318 & 0.8351 & 0.8396 & 0.8399 \\ \hline
${\bf error}_{\text{rel},1,8,4}$  & 3.137 $\times 10^{-2}$ & 1.084 $\times 10^{-2}$ & 6.899 $\times 10^{-3}$ & 1.567 $\times 10^{-3}$ & 1.234 $\times 10^{-3}$ \\ \hline \hline
$\hat S_{1,12,3}$                 & 0.8532 & 0.8748 & 0.8775 & 0.8811 & 0.8813 \\ \hline
${\bf error}_{\text{rel},1,12,3}$ & 3.291 $\times 10^{-2}$ & 8.420 $\times 10^{-3}$ & 1.290 $\times 10^{-3}$ & 1.119 $\times 10^{-3}$ & 1.067 $\times 10^{-3}$ \\ \hline
\end{tabular}
\caption{Evolution of the numerical values of the ratios $\hat S_{1,q,r}$ as defined in \eqref{Sdpq-appr} in terms of the number of iterations, for four different values of $(q,r)$.}\label{tab:S1pq}
\end{table}

\begin{figure}[!ht]
\centering
\includegraphics[width=0.5\textwidth]{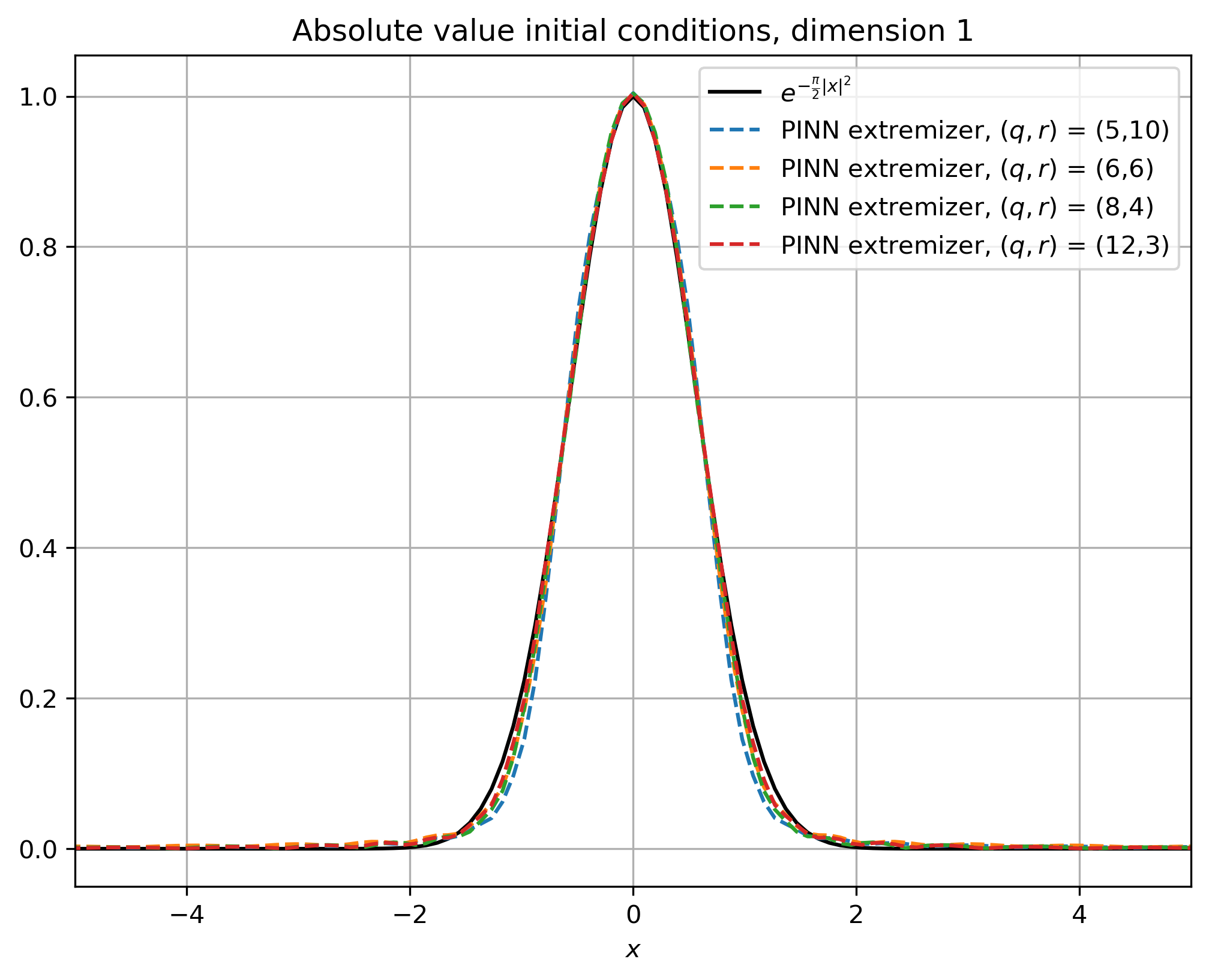}
\caption{Absolute value $|\phi_{\theta^*}(x)|$, for four different admissible pairs $(q,r)$.}
\label{fig:u0-1pq}
\end{figure}

From Table~\ref{tab:S1pq} we draw two observations. First, all experiments achieve a relative error of order $10^{-3}$, confirming that the NN algorithm recovers extremizers for the known cases $(q,r)\in\{(6,6),(8,4)\}$ and that the estimated optimal constant approaches the conjectured value \eqref{Sdpq} for the pairs $(q,r)\in\{(5,10),(12,3)\}$. Second, in almost all cases the relative error drops below $10^{-2}$ after $1000$ iterations, and the estimated constant subsequently stabilizes around $S_{1,q,r}$ (Figure~\ref{fig:S1pq}). Figure~\ref{fig:u0-1pq} compares the absolute value of the approximate extremizer with the modulus of the theoretical one (Remark~\ref{rem:unit}); in all cases, $\phi_{\theta^*}$ coincides in modulus with the gaussian $e^{-\frac{\pi}{2}|x|^2}$, so the approximate extremizers are gaussians.

\begin{figure}[t]
\centering
\begin{subfigure}[t]{0.45\linewidth}\centering
\includegraphics[width=\linewidth]{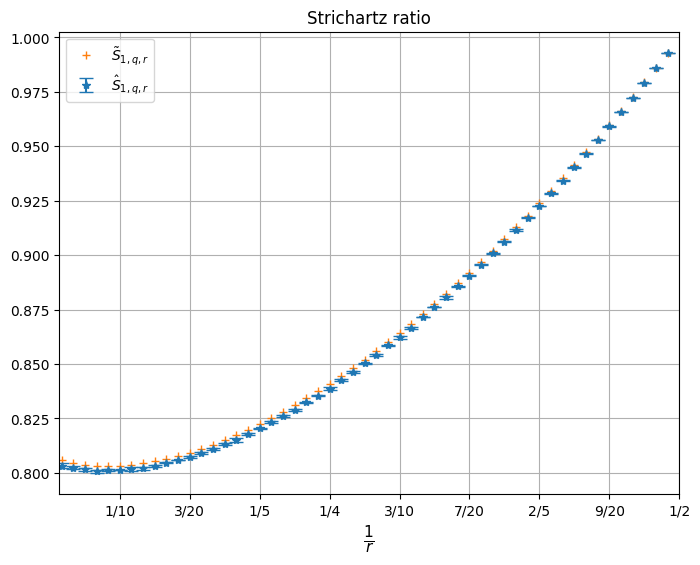}
\caption{Neral extremizer $\phi_{\theta^\star}$ ratio versus analytical gaussian ratio.}
\end{subfigure}\hfill
\begin{subfigure}[t]{0.45\linewidth}\centering
\includegraphics[width=\linewidth]{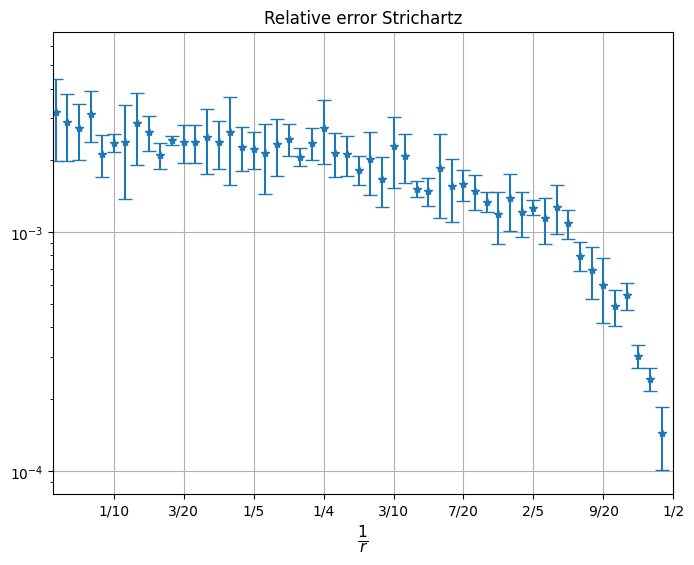}
\caption{Relative error of  $\hat{S}_{1,q,r}$, with respect to the gaussian ratio.}
\end{subfigure}
\caption{Schr\"odinger $d=1$: the neural optimization finds Gaussians as extremizers for any admissible pair.}
\label{fig:schrod-muchos-r}
\end{figure}

Additionally,  we run the full neural optimization pipeline for $59$ Schr\"odinger-admissible pairs $(q,r)$. Figure~\ref{fig:schrod-muchos-r} shows the approximate value of the Strichartz ratio via NNs, $\hat{S}_{1,q,r}$, and the relative error with respect to the Strichartz ratio of a Gaussian, $\widetilde S_{1,q,r}$. Error bars are obtained from $5$ independent runs of the algorithm for each pair and report the standard deviation. As can be seen in the figure, the neural-approximated profile corresponds to a Gaussian as well, with relative error $\lesssim 10^{-3}$ for all the admissible pairs.

\subsubsection{Dimension 1. Uneven approximations}
\begin{figure}[!ht]
\centering
\includegraphics[width=0.9\textwidth]{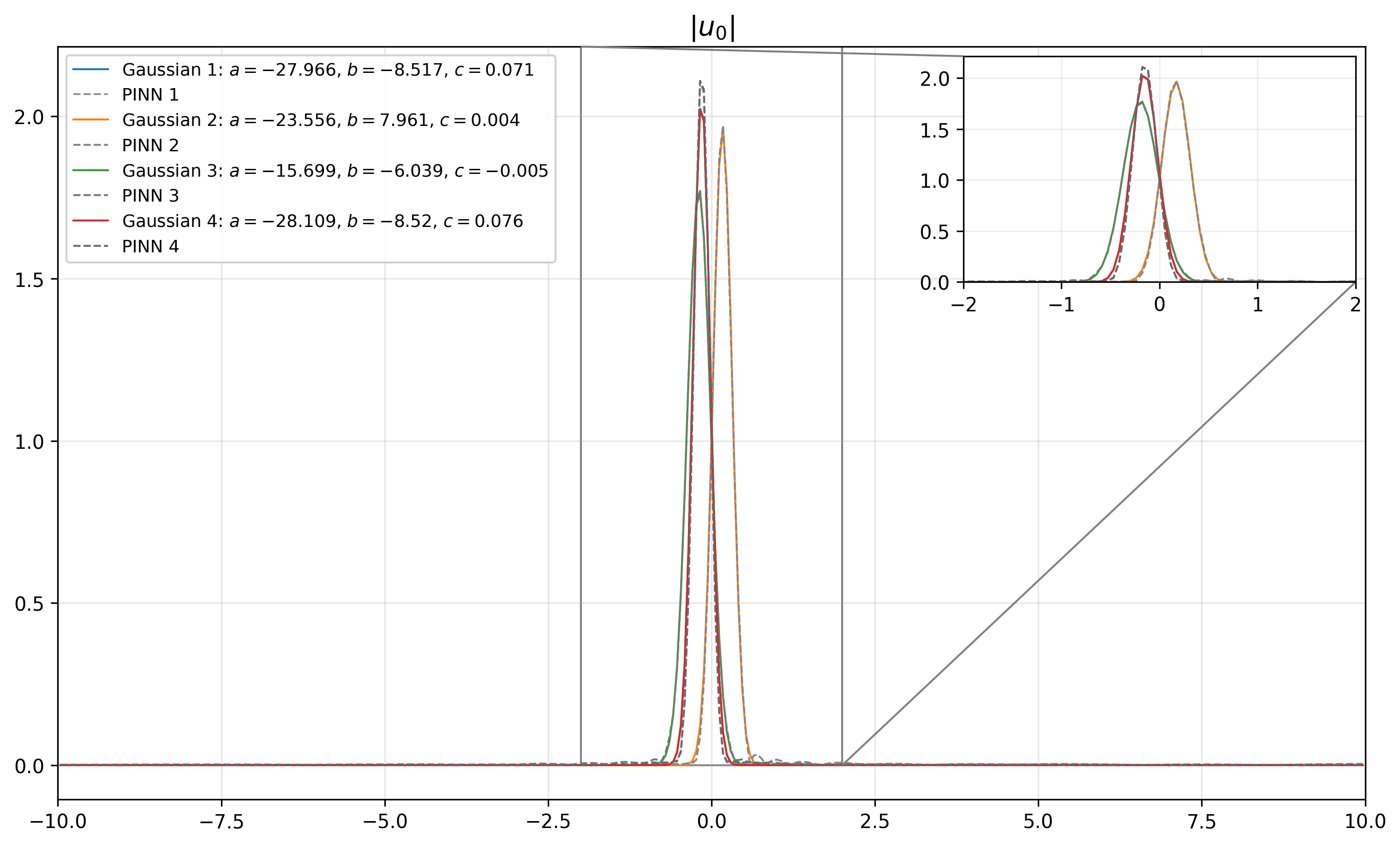}
\caption{Absolute value of uneven extremizers $|\phi_{\theta^*}(x)|$, for $(p,q)=(6,6)$ in four different realizations. For these simulations, in average $\hat{S}_{1,6,6} = 0.8128$ and ${\bf error}_{\text{rel},1,6,6} = 1.590 \times 10^{-3}$.}
\label{fig:u0-uneven}
\end{figure}

For the sake of completeness, we show that optimal neural networks remain gaussian when the parity constraint is relaxed, and that the approximate constant $\hat S_{1,q,r}$ still remains close to $S_{1,q,r}$. We only consider the case $(d,q,r)=(1,6,6)$; analogous results can be obtained for other admissible pairs. The boundary condition $|u_0(x)|=1$ is retained. Figure~\ref{fig:u0-uneven} shows the estimated extremizers in four independent realizations, each compared with a true gaussian satisfying $|u_0(x)|=1$ and $\|u_0\|_{L_x^2(\R)}=1$. Even without enforcing parity, the PINN method yields gaussians; in this case they tend to be narrower with higher amplitude.

\subsubsection{Dimension 2}
\begin{figure}[!ht]
\centering
\includegraphics[width=0.5\textwidth]{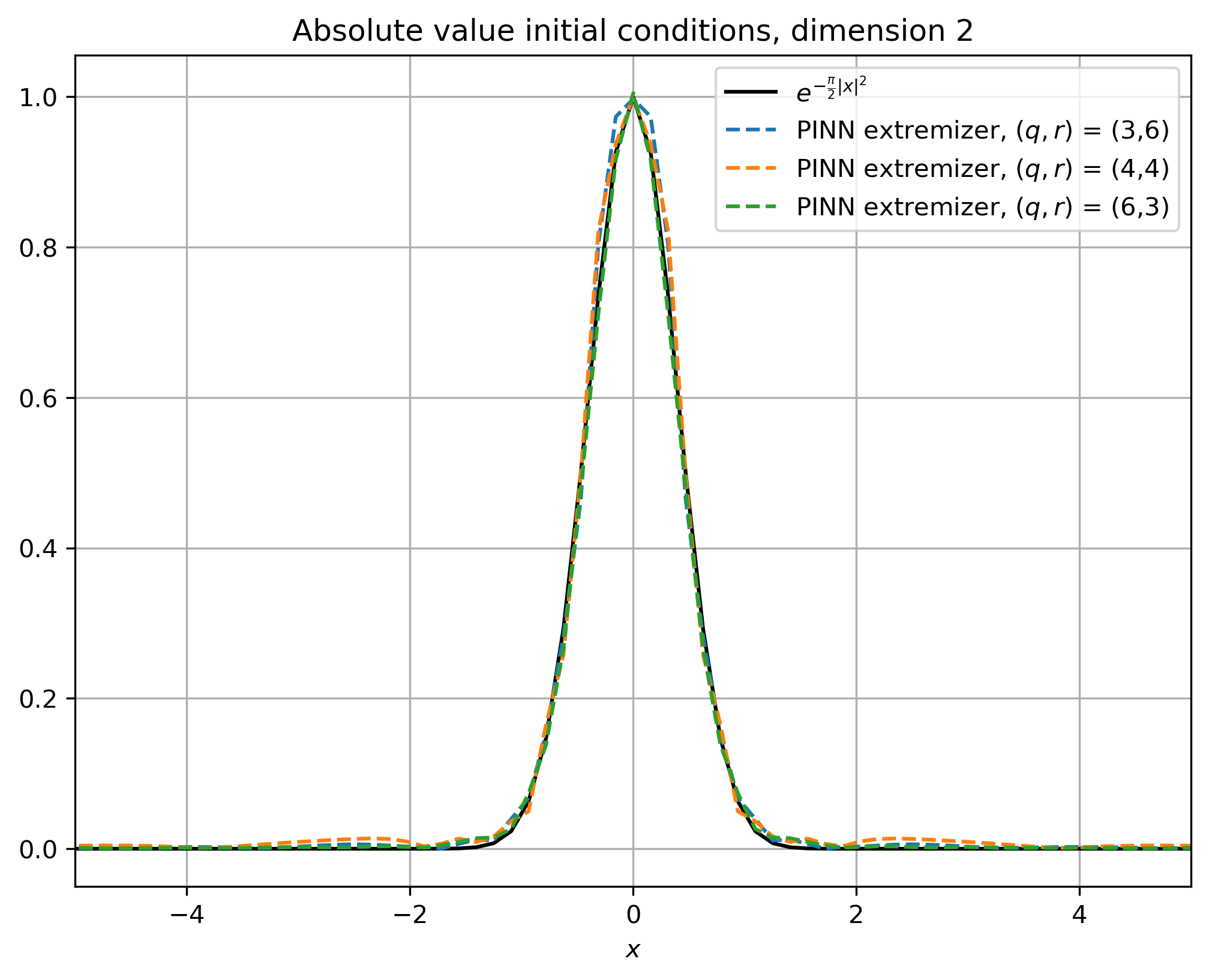}
\caption{Absolute value $|\phi_{\theta^*}(x)|$, for three different admissible pairs $(q,r)$.}
\label{fig:u0-2pq}
\end{figure}

\begin{figure}[!ht]
\centering
\begin{subfigure}[t]{0.45\linewidth}\centering
\includegraphics[width=\linewidth]{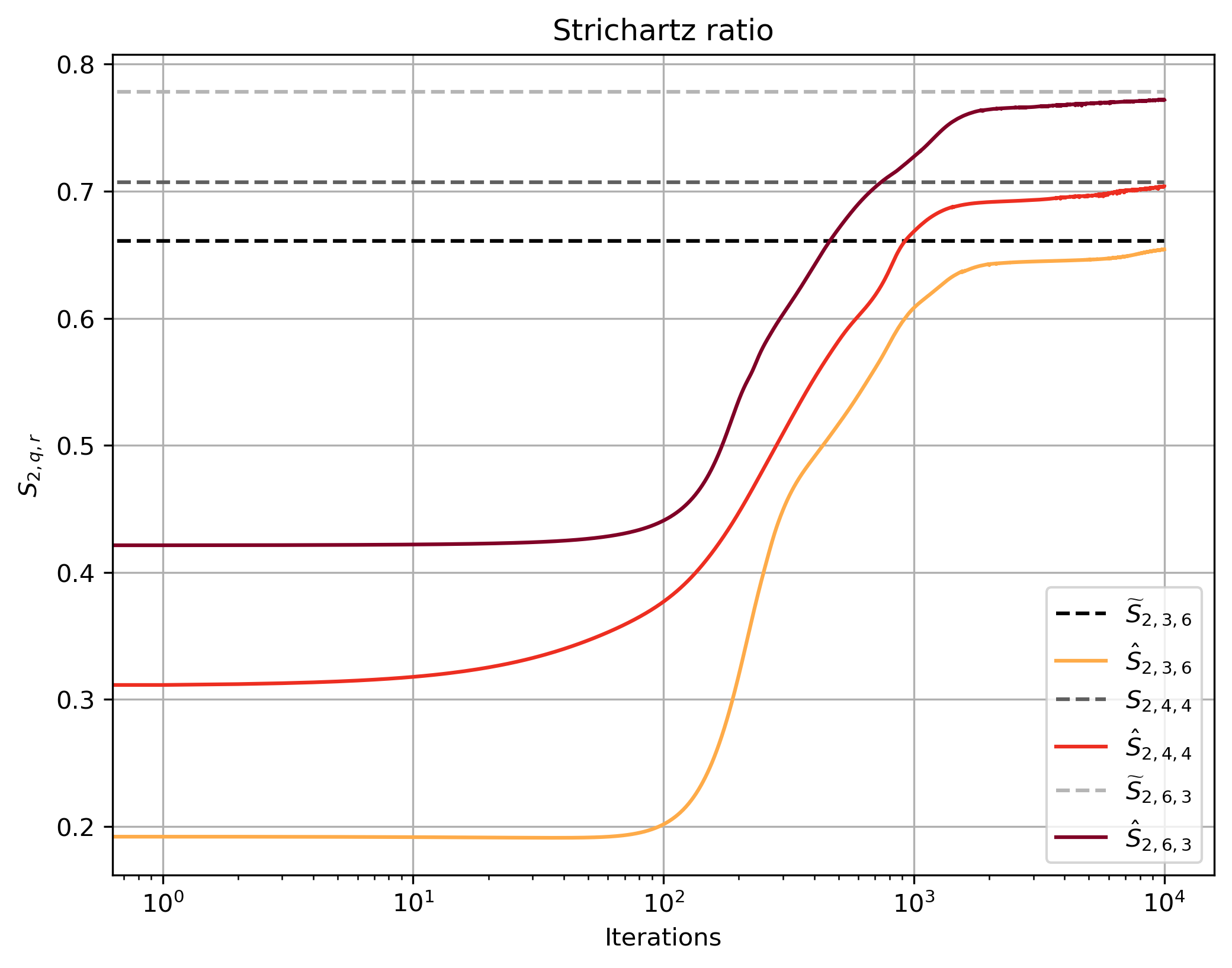}
\caption{Approximate constant $\hat{S}_{2,q,r}$.}
\end{subfigure}\hfill
\begin{subfigure}[t]{0.45\linewidth}\centering
\includegraphics[width=\linewidth]{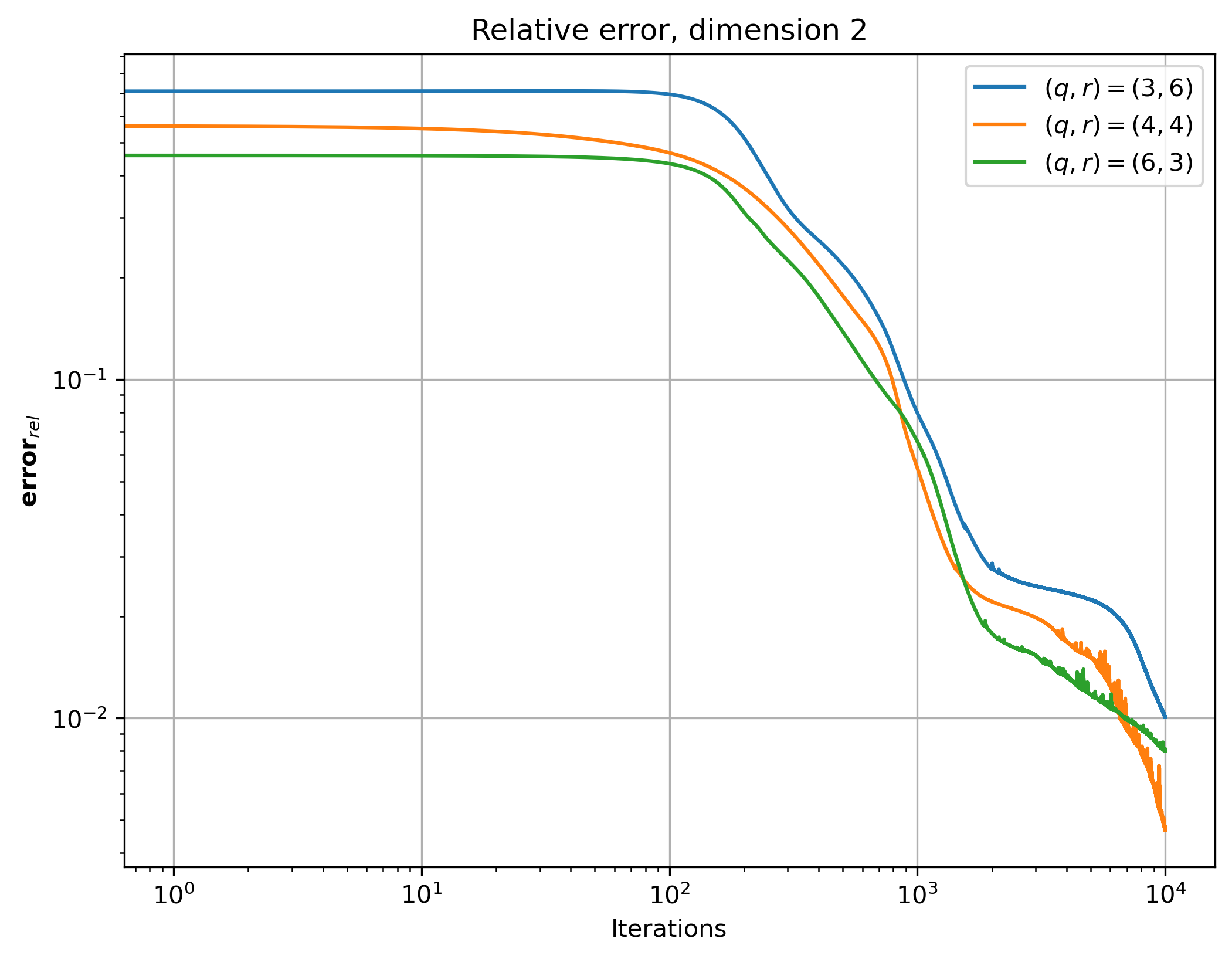}
\caption{Evolution of the relative error ${\bf error}_{\text{rel},2,q,r}$.}
\end{subfigure}
\caption{Schr\"odinger $d=2$: NN results for three different admissible pairs $(q,r)$.}
\label{fig:S2pq}
\end{figure}

For this dimension, the functionals $\mathcal J_{L_x^2(\R^d)}$ and $\mathcal J_{L_t^qL_x^r(\R^{d+1})}$ are computed with $R = 10$, $T=1$, $M=32$ and $N=128$. Due to the higher computational cost in dimension 2, the number of grid points was reduced. We first consider the pair with known extremizer, $(d,q,r) = (2,4,4)$. Additionally, we consider two further experiments with admissible pairs $(d,q,r) \in \{(2,3,6),(2,6,3)\}$. The results are summarized in Table \ref{tab:S2pq} and Figures \ref{fig:u0-2pq}, and \ref{fig:S2pq}. As in dimension 1, the values in Table \ref{tab:S2pq} are averages of 5 independent realizations of the algorithm. From Table \ref{tab:S2pq} we can see that $\hat{S}_{2,q,r}$ is close to $S_{2,q,r}$ given in \eqref{Sdpq} (relative error below $10^{-2}$), but in contrast to dimension~1, convergence here is slower. Indeed, for $d=2$, the relative error is still above $10^{-2}$ after 5000 iterations. This can be seen in Figure \ref{fig:S2pq} as well.

Figure \ref{fig:u0-2pq} reports the value of $|\phi_{\theta^*}|$. Since $\phi_{\theta^*}$ has inputs in $\R^2$, Figure \ref{fig:u0-2pq} shows the absolute value projected on the diagonal $x_1 = x_2$. As shown in the figure, the discovered extremizers are Gaussian, and comparable with $e^{-\frac\pi2 |x|^2}$. However, in this dimension the approximation of the Gaussian's decaying tail is worse than in dimension $d=1$. It is important to note that the grid spacing $\Delta x$ considered in dimension 1 is smaller than the one considered in dimension 2. Despite being less accurate than in dimension 1, the results are still reasonable, with relative error below $10^{-2}$, demonstrating the capability of neural networks for approximating the optimal constant $S_{2,q,r}$ for Schrödinger-admissible pairs.

\begin{table}[!ht]
\def\arraystretch{1.4}
\begin{tabular}{|c||c|c|c|c|c|}
\hline
Iterations                        & 100 & 500 & 1000 & 5000 & 10000 \\ \hline \hline
$\hat S_{2,3,6}$                 & 0.2322 & 0.6093 & 0.6391 & 0.6535 & 0.6621 \\ \hline
${\bf error}_{\text{rel},2,3,6}$ & 6.487 $\times 10^{-1}$ & 7.804 $\times 10^{-2}$ & 3.303 $\times 10^{-2}$ & 1.123 $\times 10^{-2}$ & 6.361 $\times 10^{-3}$ \\ \hline \hline
$\hat S_{2,4,4}$                  & 0.3651 & 0.6617 & 0.6855 & 0.6981 & 0.7067 \\ \hline
${\bf error}_{\text{rel},2,4,4}$  & 4.836 $\times 10^{-1}$ & 6.424 $\times 10^{-2}$ & 3.058 $\times 10^{-2}$ & 1.268 $\times 10^{-2}$ & 6.949 $\times 10^{-3}$ \\ \hline \hline
$\hat S_{2,6,3}$                  & 0.4450 & 0.7355 & 0.7602 & 0.7692 & 0.7727 \\ \hline
${\bf error}_{\text{rel},2,6,3}$  & 4.282 $\times 10^{-1}$ & 5.496 $\times 10^{-2}$ & 2.327 $\times 10^{-2}$ & 1.163 $\times 10^{-2}$ & 7.100 $\times 10^{-3}$ \\ \hline
\end{tabular}
\caption{Evolution of the numerical values of the ratios $\hat S_{2,q,r}$ as defined in \eqref{Sdpq-appr} in terms of the number of iterations, for three different values of $(q,r)$.}\label{tab:S2pq}
\end{table}


{
From the findings regarding the Strichartz estimates in the case of the linear Schrödinger model, we find supporting evidence for the validity of the well-known conjecture on the extremizers of \eqref{Str-Sch}:    In dimension $d =1,2$, Strichartz estimates with admissible pair $(p,q)$ are extremized by Gaussians. Moreover, the optimal constant is
    \be\label{Sdpq}
    S_{d,q,r} = \widetilde{S}_{d,q,r},
    \ee
with $\widetilde{S}_{d,q,r}$ given in \eqref{tildeS}.

\begin{remark}
    Note that for $(d,q,r) \in \{(1,6,6),(2,4,4),(1,8,4)\}$, the exact value of the optimal constant, given in Table \ref{tab:Strichartz-ctes}, coincides with the provided formula in the previous conjecture.
\end{remark}
}

\section{Airy models}
\label{app:estimates}

The Airy initial-value problem
\begin{equation}\label{Airy}
\partial_t v + \partial_x^3 v = 0,\qquad v(0,x)=v_0(x),
\end{equation}
is the linear part of the (focusing) generalized Korteweg--de Vries family $\partial_t u+\partial_x(\partial_x^2 u+u^p)=0$, $p\ge 2$ \cite{KdV,KPV1,Bourgain,CKSTT,Guo,Kishimoto,KKSV,KiVi19}, whose soliton solutions are
\begin{equation}\label{soliton}
Q_{c,p,x_0}(t,x)\;:=\;Q_{p,c}(x-ct-x_0),\qquad Q_{p,c}(s)=c^{1/(p-1)}Q_p(\sqrt{c}s),
\end{equation}
with $Q_p$ as in \eqref{Qk}; multi-soliton solutions for $p=2,3$ are constructed in \cite{Hir71,Hir72,Martel}. The completely integrable mKdV ($p=3$) admits the explicit four-parameter breather \cite{Wadati,Lamb,AM}
\begin{equation}\label{BBB}
B(t,x;\alpha,\beta,x_1,x_2)\;=\;2\sqrt{2}\,\partial_x\arctan\!\left(\frac{\beta}{\alpha}\,\frac{\sin(\alpha x+\delta t+x_1)}{\cosh(\beta x+\gamma t+x_2)}\right),
\end{equation}
where $\gamma=\beta(3\alpha^2-\beta^2)$, $\delta=\alpha(\alpha^2-3\beta^2)$. The full literature is reviewed in Section~\ref{sec:prelim}; see also \cite{KPV2,Ale1,GoPe,Na1,Na2,AGV,AM,Munoz,AleProy1,AleProy2,OW,Ale,Ale0,ChenLiu,Munoz2,AFM,AM_book,Ka,Bous,FPU} for breather and gKdV background.

\subsection{Airy-Strichartz inequalities}
 The Airy--Strichartz framework \eqref{AS}--\eqref{Aqr} is recalled in Section~\ref{sec:intro}; we expand here on the explicit ratio computations for solitons and breathers. We will work then in a more general framework. Indeed, the following Airy-Strichartz estimates are well-known:
\be\label{airy-strichartz}
\| |D_x|^\gamma (e^{-t\partial_x^3}u)\|_{L^q_tL^r_x(\R\times\R)} \leq C \|u\|_{L^2_x},
\ee
where
\[
2 \leq q,r< \infty, \qquad -\gamma +\frac 3q + \frac 1r = \frac12, \qquad -\frac 12 < \gamma \leq \frac 1q.
\]
Here, $|D_x|^\gamma=\mathcal F^{-1} (|\xi|^\gamma  \mathcal F (\cdot))$ whenever it is well-defined. For the Airy-Strichartz estimates \eqref{airy-strichartz}, following the literature, we define the sharp constant as
\be\label{Aqr-app}
A_{q,r} := \sup_{u_0 \in L_x^2(\R)-\{0\}} \frac{\||D_x|^{\gamma}(e^{-t\partial_x^3}u_0)\|_{L_t^qL_x^r(\R_t\times\R_x)}}{\|u_0\|_{L_x^2(\R_x)}},
\ee
provided it exists. The existence of extremizers for Airy-Strichartz estimates is proved in \cite{FS18} for the critical case $\gamma = \frac 1q$, however without any insights about the shape of the extremizers. In such a critical case, however, an important result can be stated regarding extremizers itself, and a lower bound for $A_{q,r}$ from $S_{1,q,r}$, as defined in \eqref{Sdqp}. Adapted to our notation, we have the following
\begin{theorem}[\cite{FS18}]
Let $4 < q < \infty$ and $r$ such that $(1,q,r)$ is Schrödinger admissible. Then, all maximizing sequences for $A_{q,r}$ are precompact up to symmetries if and only if
\be\label{lb-AS}
A_{q,r} > 3^{-\frac1q}a_r S_{1,q,r} =: \widetilde A_{q,r},
\ee
where $A_{q,r}$ is defined in \eqref{Aqr}, $S_{1,q,r}$ is defined in \eqref{Sdqp}, and
\[
a_r := \frac{2^{1/2}}{\pi^{1/(2r)}} \left(\frac{\Gamma\left(\frac{r+1}{2}\right)}{\Gamma\left(\frac{r+2}{2}\right)}\right)^{1/r} > 1.
\]
In particular, if \eqref{lb-AS} holds, then there is a maximizer of $A_{q,r}$.
\end{theorem}

\begin{figure}[!ht]
\centering
\includegraphics[width=0.5\textwidth]{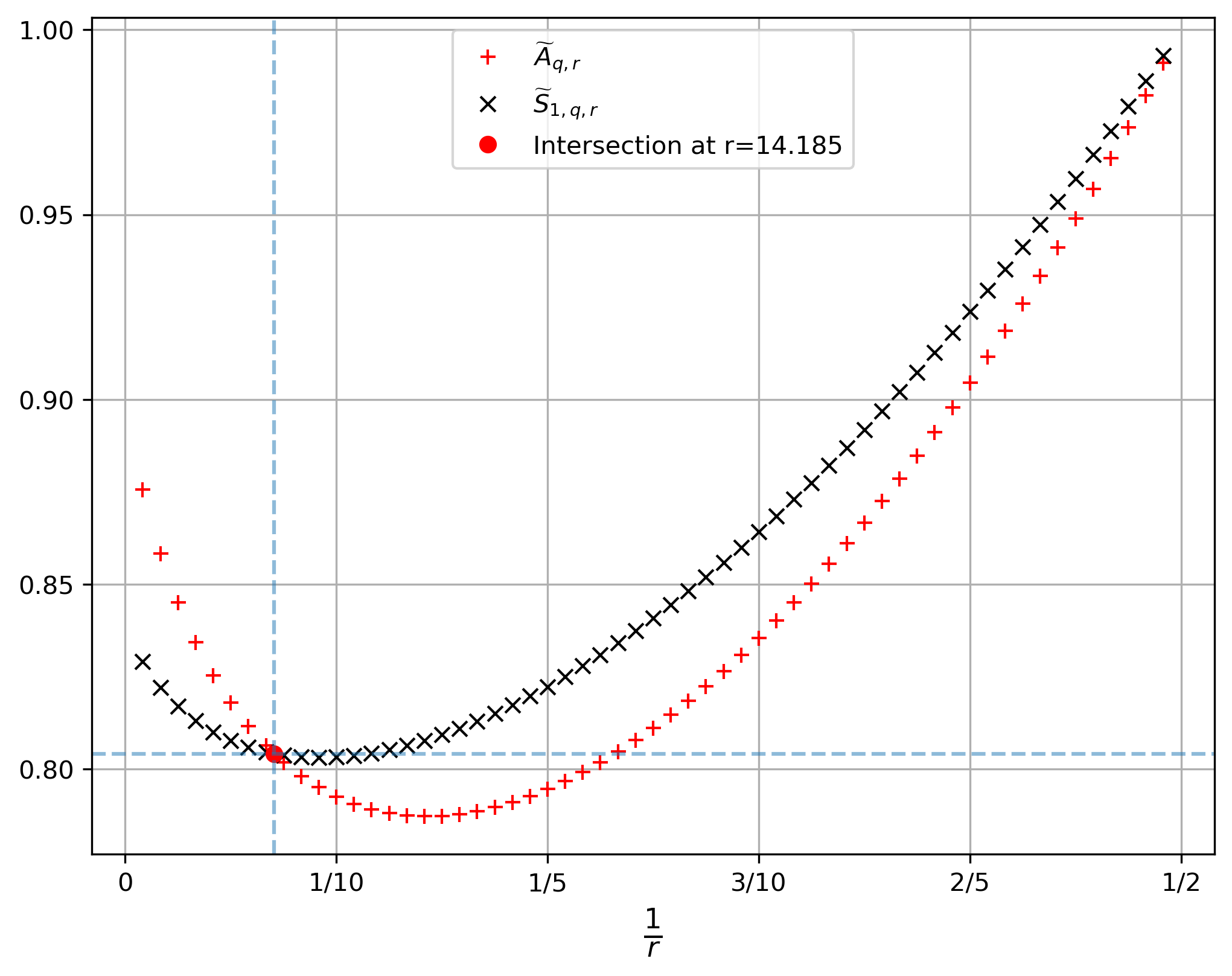}
\caption{Comparison of Airy-Strichartz (Airy) and Strichartz (Schrödinger) sharp constants.}
\label{fig:comp-AS-S}
\end{figure}

\begin{figure}[ht!]
    \centering
    \includegraphics[width=0.7\linewidth]{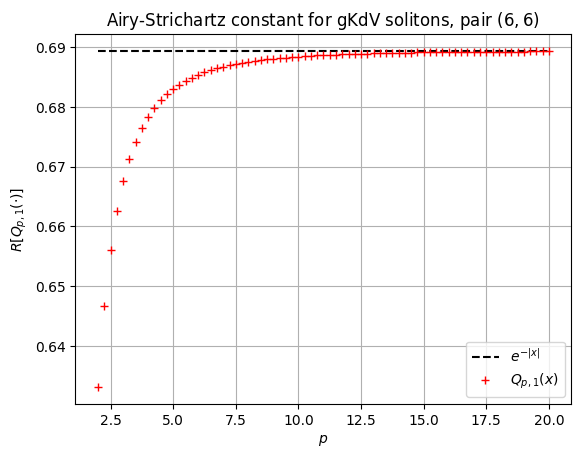}
    \caption{In red: Airy-Strichartz ratio of initial conditions $Q_{p,1}(x)$ for different values of $p$. In black, the ratio for the formal scaling-limit profile $e^{-|x|}$.}
    \label{fig:ratio-AS-sol}
\end{figure}

\begin{figure}[ht!]
    \centering
    \includegraphics[width=0.7\linewidth]{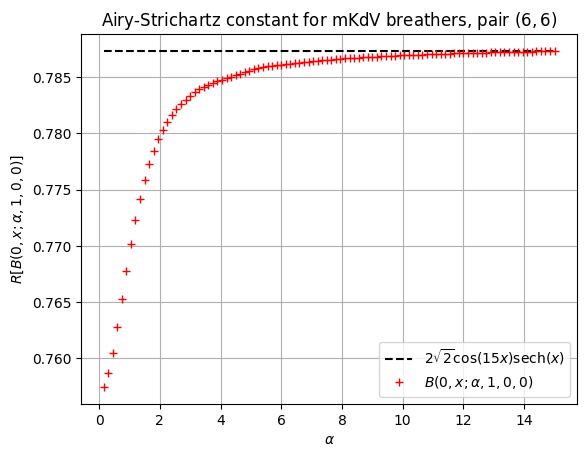}
    \caption{In red: Airy-Strichartz ratio of initial conditions $B(0,x;\alpha,1,0,0)$ for different values of $\alpha$. In black, the ratio for the limiting breather $2\sqrt{2} \cos(\alpha x)\sech(x)$.}
    \label{fig:ratio-AS-breather}
\end{figure}

\begin{figure}[ht!]
    \centering
    \includegraphics[width=0.7\linewidth]{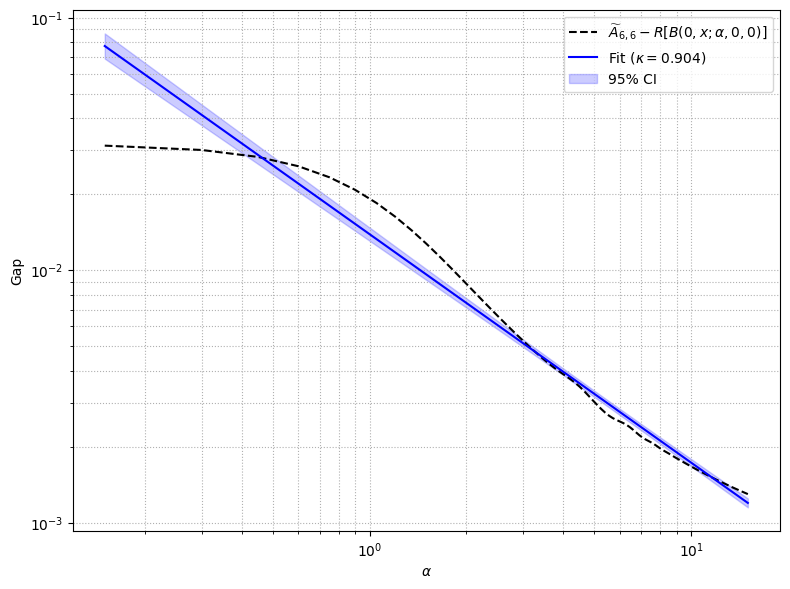}
    \caption{In black, difference between $\widetilde A_{6,6}$ and $R[B(0,x;\alpha,1,0,0)]$, for increasing $\alpha$. In blue, power-law decay $\sim 0.0139\,\alpha^{-0.9041}$.}
    \label{fig:fit-slope-breather}
\end{figure}

One can easily compute the value for $\widetilde A_{q,r}$, corresponding to a hypothetical lower bound for $A_{q,r}$, given by the inequality \eqref{lb-AS}. Figure \ref{fig:comp-AS-S} summarizes the behavior of the constants $\widetilde S_{1,q,r}$ and $\widetilde A_{q,r}$ as a function of $\frac1r$. As can be deduced from this figure, for small values of $r$, the sharp constant $\widetilde A_{q,r}$ is below $\widetilde S_{1,q,r}$, while it becomes larger for $r > 14.185$ approximately. Let us note that in the known particular case $(q,r)=(6,6)$, condition \eqref{lb-AS} should read
\be\label{A66}
A_{6,6} > \widetilde A_{6,6} \approx 0.7886,
\ee
revealing that $A_{6,6}$ must be strictly greater than this value for an $L^2$ extremizer to exist. When $A_{q,r}$ collapses to the universal lower bound $\widetilde A_{q,r}$, the precompactness criterion of \cite{FS18} fails and no $L^2$ extremizer can be attained.
To have a notion of the validity of this condition, an interesting fact is to compare it with the behavior of solitons \eqref{soliton} in \eqref{airy-strichartz}. One can easily check that the ratio
\be\label{ratio_c}
R[Q_{p,c}(\cdot -x_0)]:=\frac{\| |D_x|^\gamma (e^{-t\partial_x^3}Q_{p,c})\|_{L^q_tL^r_x}}{\|Q_{p,c}\|_{L^2_x}}
\ee
is independent of the scaling $c$ and shift $x_0$. This is a consequence of the fact that if $u(t,x)$ solves Airy, then for any $\lambda$ and $c>0$, $\lambda u(c^{\frac32}t,c^{\frac12}x)$ also does. Therefore, the independence of the scaling property is more general than just \eqref{ratio_c}. In terms of numerical computations, this will allow us to discard huge regions of parameters by only considering certain regions. The ratio of $Q_{p,1}(\cdot)$ for different values of $p$ is given in Figure \ref{fig:ratio-AS-sol}. In addition, we compare them with the ratio associated with the infinite scaling limit soliton $R[e^{-|\cdot|}]$. This function is obtained by formally taking limits to infinity on the scaling associated with each soliton. We realize that the ratio of the soliton is increasing with $p$, does not exceed the ratio of $e^{-|\cdot|}$, and numerically converges to the limiting ratio as $p$ increases. We thus have a first lower bound on $A_{q,r}$, which in the $(q,r)=(6,6)$ case is below the expected value $\widetilde A_{q,r}$ in \eqref{A66}. 
In the case of breathers \eqref{BBB}, following the same previous scaling and shift considerations, we shall have from \eqref{BBB} that
\be\label{ratio_B}
\begin{aligned}
R[B(0,\cdot;\alpha,\beta,x_1,x_2)] = &~{} \frac{1}{2\sqrt{2\beta}} \left\| |D_x|^\gamma (e^{-t\partial_x^3}B(0,\cdot ;\alpha,\beta,x_1,x_2)) \right\|_{L^q_tL^r_x}\\
= &~{} \frac{1}{2\sqrt{2}} \left\| |D_x|^\gamma \left( e^{-t\partial_x^3}B\left(0,\cdot ;\frac{\alpha}{\beta},1,0, 0 \right)\right) \right\|_{L^q_tL^r_x}.
\end{aligned}
\ee
Here we have used that the (conserved) mass $\frac12\int B^2(t,x)dx$ of the breather is $4\beta$ \cite{AM,AM_book}, and the identity
\[
B(0,x;\alpha,\beta,x_1,x_2) = \beta B\left(0,\beta x;\frac{\alpha}{\beta},1,x_1,x_2 \right),
\]
along with $L^q_tL^r_x$ invariance under shifts in time and space. Therefore, \eqref{ratio_B} tells us that, in the case of breathers, we can assume that the $\beta$ scaling is equal to one, and the $\alpha$ scaling is a positive free parameter. The ratio of $B(0,\cdot;\alpha,1,0,0)$ for different values of $\alpha$ is summarized in Figure \ref{fig:ratio-AS-breather} for the pair $(6,6)$. We deduce that the ratio of the breather is increasing with $\alpha$. Notice that when $\alpha \gg 1$ the limit of the breather behaves as
\[
B(0,x;\alpha,1,0,0) \approx 2\sqrt{2}\cos(\alpha x) \sech(x) + O\left(\frac 1\alpha\right).
\]
Figure \ref{fig:ratio-AS-breather} shows in black the behavior of the limiting breather $2\sqrt{2}\cos(\alpha x) \sech(x)$ when $\alpha = 15$. Note from the figure that the ratios for breathers are below those of the limiting breather, and all considered ratios are below the lower bound \eqref{lb-AS}-\eqref{A66}, but increasingly close as $\alpha$ grows to infinity. It is important to notice that the Airy evolution and the nature of breathers are highly oscillatory, and it is difficult to capture all the effective frequencies when performing FFT. Later, Figure \ref{fig:fit-slope-breather} shows that the breather ratio approaches $\widetilde A_{6,6}$ as $\alpha \to \infty$. A log-log least-squares fit reveals a power-law decay of the gap of the form $C\alpha^{-\kappa}$, with exponent $\kappa \approx 0.904$. The fit is performed in the asymptotic regime and is consistent with a monotone collapse of the gap.

\subsection{The case of Hermite functions}
\label{app:hermite}
{To establish an independent analysis, we focus on the Schrödinger admissible pair $(q,r)=(6,6)$ in the critical case $\gamma=\frac1q$. In order to complement the study of ratios for the Airy-Strichartz inequality, we consider the orthonormal $L^2$ basis of Hermite functions. The Hermite functions are formally defined for $n \in \N$ as
\[
f_n(x) = He_n(x) \frac{e^{-\frac{x^2}{4}}}{\sqrt{n!}(2\pi)^{\frac14}},
\]
where $He_n(\cdot)$ is the probabilistic Hermite polynomial
\[
He_n(x) = (-1)^n e^{\frac{x^2}{2}} \frac{d^n}{dx^n} e^{-\frac{x^2}{2}}.
\]
We consider the first 5 Hermite functions and we compute the ratios for those functions as initial conditions, by changing the value of the time-domain boundaries. We consider a grid of 30 uniform times $T_\ell$ between 0.1 and 15. Then we compute the Airy-Strichartz norm on the domain $[-T_\ell,T_\ell]\times[-R,R]$, with $R\gg1$ fixed. In our computations, $R$ is chosen as $R=2500$ to avoid periodicity issues on the FFT. Figure \ref{fig:AS-Hermite-time-exp} outlines the ratio of the Hermite functions when the time domain is increasing. As shown in the figure, for the five studied functions, their ratios are stable at time $T_\ell = 15$. Using this time as boundary, we can quantify the Airy-Strichartz ratio for a bigger family of Hermite functions. In particular, Figure \ref{fig:AS-Hermite} shows the ratios when the initial condition is chosen as $f_n(x)$, for $n \in \{0,1,\ldots,20\}$. From the figures we can deduce that the Hermite function $f_1$ is the Hermite function with bigger Airy-Strichartz ratio, but it is still far away from the lower bound $\widetilde{A}_{6,6}$.

\begin{figure}[ht!]
    \centering
    \includegraphics[width=0.5\linewidth]{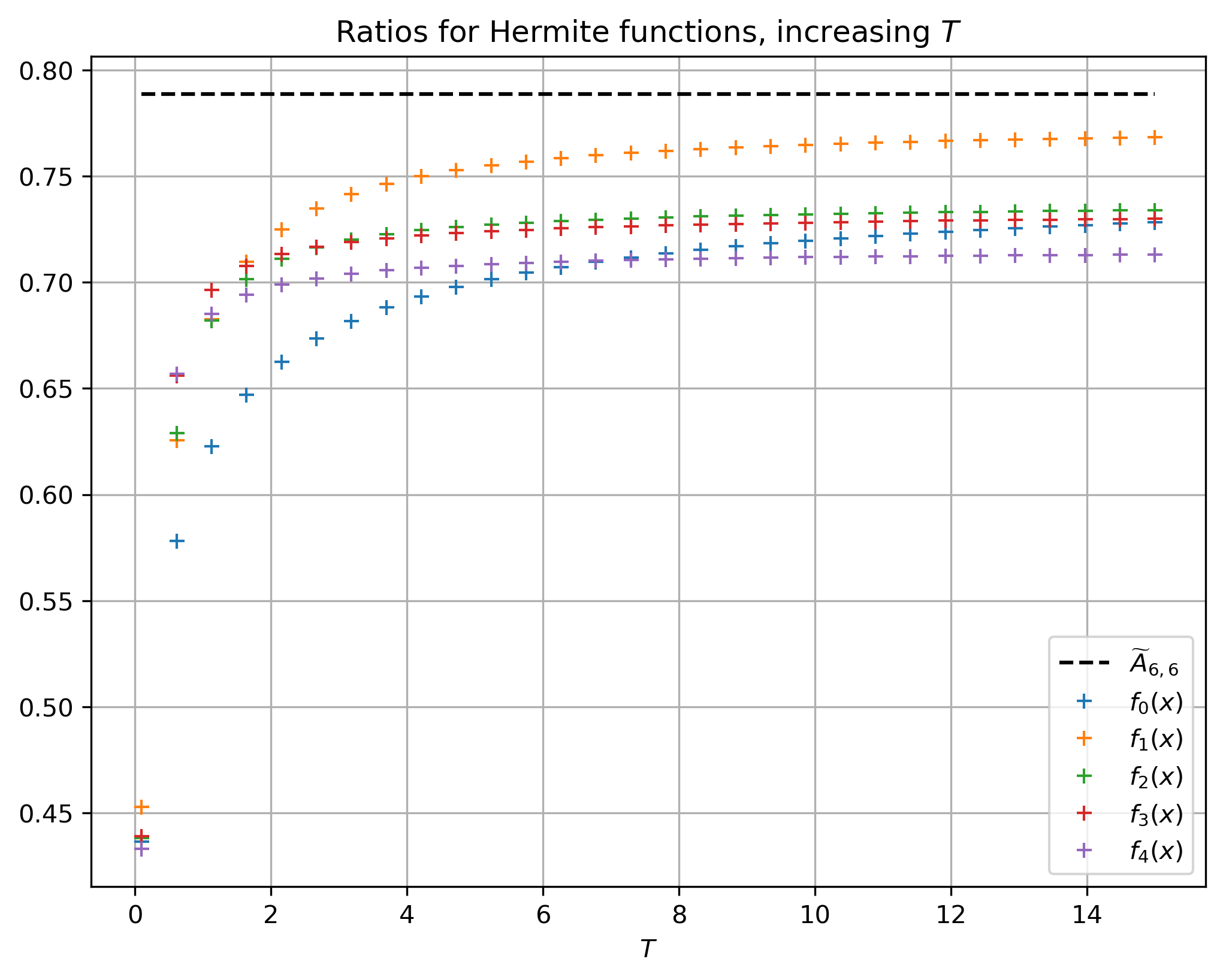}
    \caption{In black, the lower bound $\widetilde{A}_{6,6}$. In ``plus" dots, the ratios for Hermite functions with $n=0,1,2,3,4$, when increasing time boundary $[-T,T]$x}
    \label{fig:AS-Hermite-time-exp}
\end{figure}

\begin{figure}[ht!]
    \centering
    \includegraphics[width=0.5\linewidth]{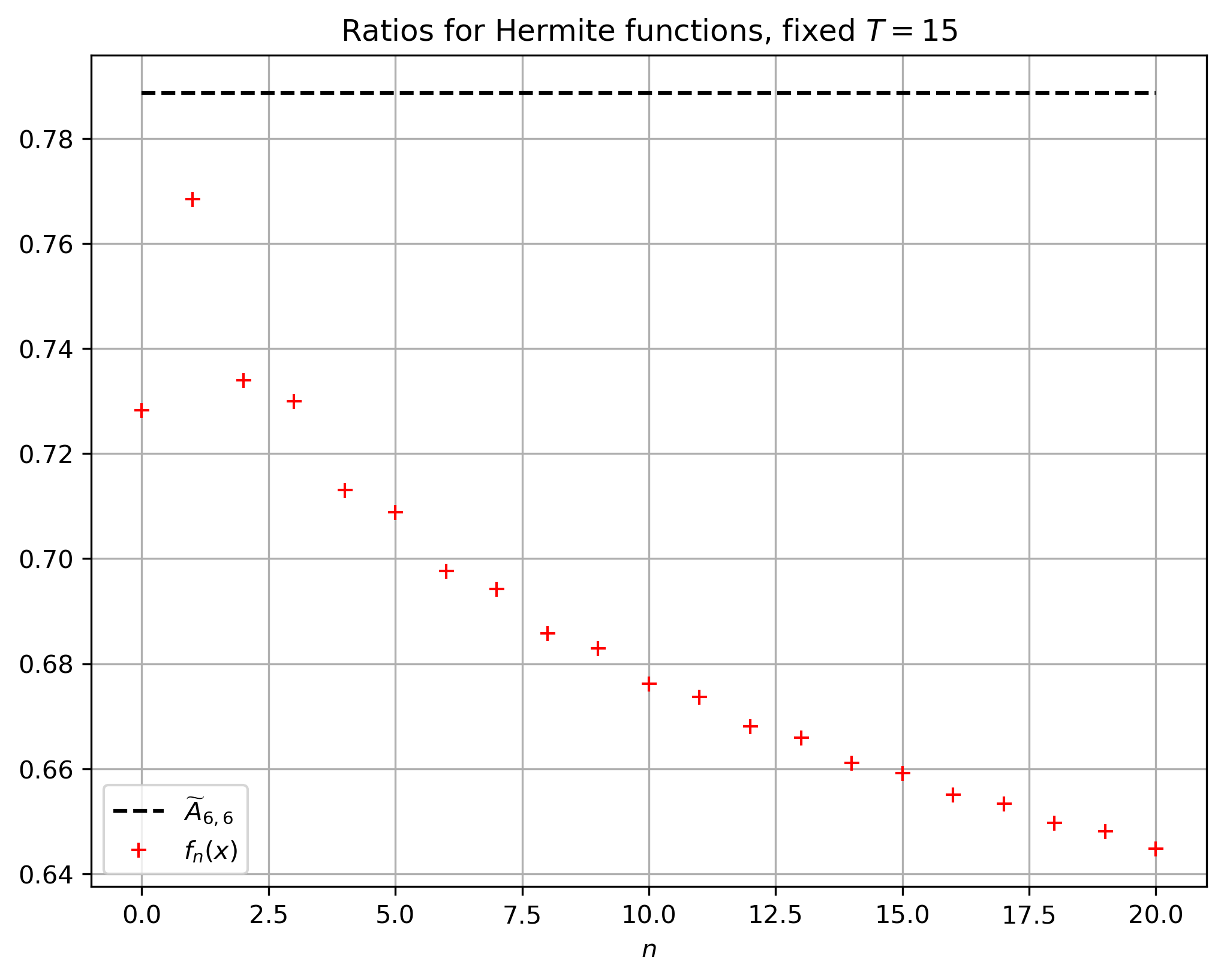}
    \caption{In black, the lower bound $\widetilde{A}_{6,6}$. In red, the ratios for Hermite functions with $n \in \{0,1,\ldots,20\}$.}
    \label{fig:AS-Hermite}
\end{figure}

Before attacking the problem with neural networks, we use the $L_x^2$ Hermite basis to propose an Ansatz of the extremizer having the form
\be\label{linear-comb}
u_0(x) = \sum_{n=0}^{20} b_n f_n(x),
\ee
with $b_n \in \R$ being free parameters. We have chosen 20 functions to ensure sufficient expressivity of the obtained $u_0$, but any other similar number will also be enough. We employ an optimization algorithm to obtain the optimal values $b_n^*$ that maximize the Airy Strichartz ratio; in other words, we seek the optimal sequence $(b_n^*)_{n=0}^{20}$ such that
\[
(b_n^*)_{n=0}^{20}=\arg\min_{(b_n)_{n=0}^{20}}- \frac{\||D^{\frac16}_x|u_0\|_{L_t^6L_x^6(\R\times\R)}}{\mathcal \|u_0\|_{L_x^2(\R)}}.
\]
By performing the optimization algorithm, we have obtained the following facts: The numerically found coefficients are shown in Table \ref{tab:bn}. By plugging these coefficients into $u_0$ in \eqref{linear-comb}, we obtain an initial condition that is comparable with a breather, as presented in Figure \ref{fig:hermite-linear-comb}, and furthermore, the ratio obtained by the optimization algorithm corresponds to 
\[
\hat A_{6,6}^{\text{Hermite}} = 0.7828 < \widetilde A_{6,6}.
\]
It is also noted that after the value $n=13$, the coefficients found decrease and have values that are reasonably and increasingly smaller than the previous ones. 
\begin{table}[!ht]
\centering
\begin{tabular}{|c|c||c|c|}
\hline
$n$ & $b^*_n$ & $n$ & $b^*_n$ \\
\hline
0  & -0.0327 & 11 & -4.1902 \\ \hline
1  & 0.7151  & 12 & 0.4417 \\ \hline
2  & 0.3791  & 13 & -3.5822 \\ \hline
3  & -1.8792 & 14 & -0.3319 \\ \hline
4  & -1.2148 & 15 & 0.2587 \\ \hline
5  & 1.2915  & 16 & 0.7547 \\ \hline
6  & -0.7746 & 17 & 0.9837 \\ \hline
7  & 5.0181  & 18 & 0.9870 \\ \hline
8  & 3.0026  & 19 & 1.0004 \\ \hline
9  & 3.0744  & 20 & 0.9995 \\ \hline
10 & 3.0210  &    &        \\
\hline
\end{tabular}
\caption{Coefficients $b^*_n$ for the function $u_0$ defined in \eqref{linear-comb}.}
\label{tab:bn}
\end{table}

\begin{figure}[ht!]
    \centering
    \includegraphics[width=0.7\linewidth]{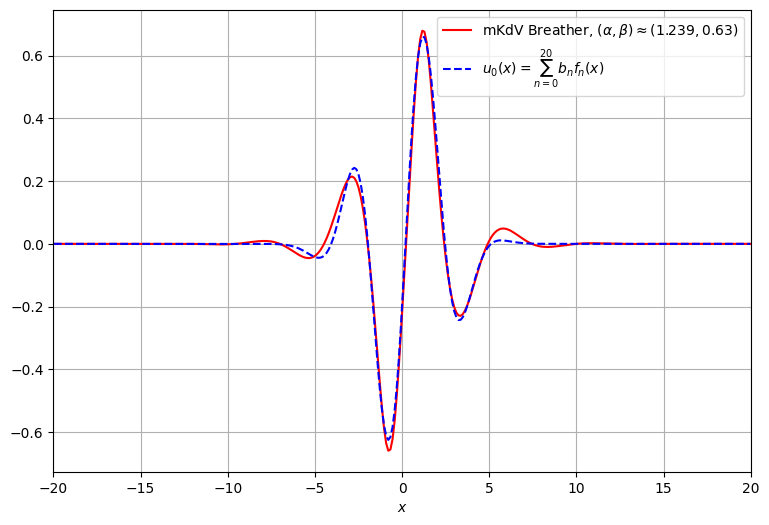}
    \caption{In blue, the function $u_0(x)$ defined in \eqref{linear-comb} with coefficients $b_n^*$ given in Table \ref{tab:bn}. In red, mKdV breather with $\alpha \approx 1.239$, $\beta \approx 0.630$, $t_0 \approx 0.313$, $x_0 \approx -1.591$.}
    \label{fig:hermite-linear-comb}
\end{figure}
}
\subsection{Numerical experiments}\label{app:AS-numerics}
{For the Airy model, we consider the critical case $\gamma = \frac 1q$, with the same Schrödinger-admissible pairs considered in the Schrödinger model. All experiments use 4 hidden layers and 20 neurons per layer. As a new insight, the NN algorithm produces an approximate extremizer whose structure is closely related to an mKdV breather, as can be seen in the left panel of Figure \ref{fig:u0-kdv}. The right panel of the same figure summarizes the evolution of the constant $\hat A_{q,r}$ compared with the lower bound $\widetilde A_{q,r}$ defined in \eqref{lb-AS}, across iterations. From the figure, one can deduce the following: first, the width of the optimal neural network is highly related to the exponent $q$; and second, the approximated constant $\hat{A}_{q,r}$ stabilizes near the lower bound $\widetilde{A}_{q,r}$.

\begin{figure}[!ht]
\centering
\includegraphics[width=0.45\textwidth]{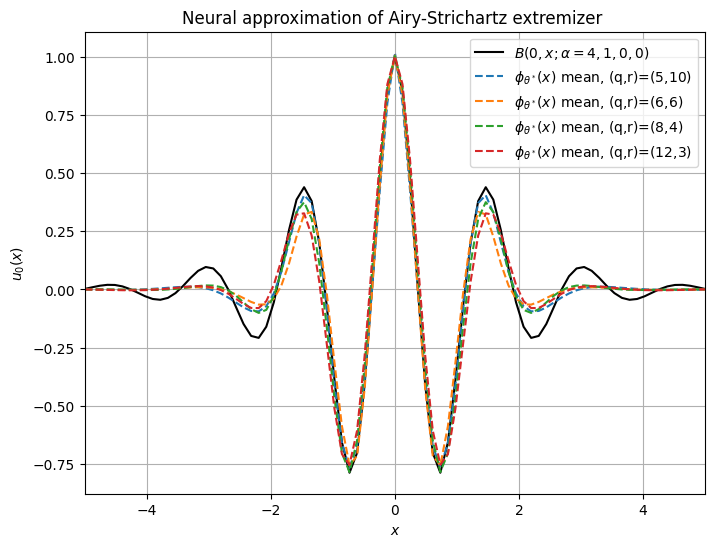}
\includegraphics[width=0.45\textwidth]{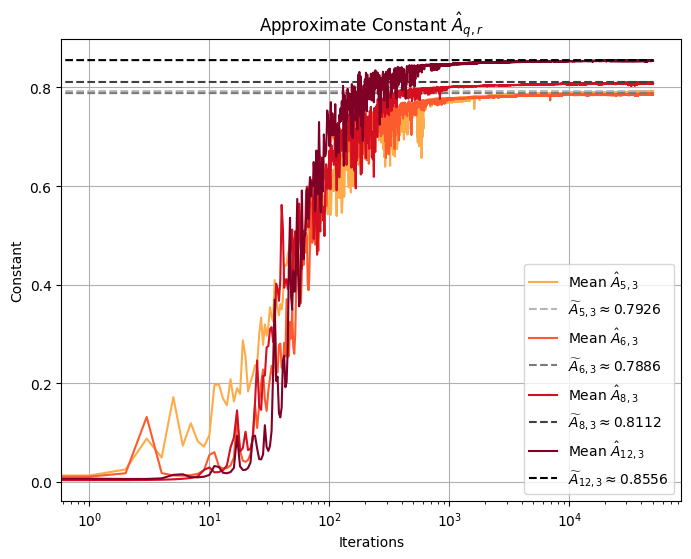}
\caption{Left: Neural extremizer obtained by means of $5$ realizations, for each pair $(q,r)$, compared with the mKdV breather $B(0,x;4,1,0,0)$. Right: Evolution of constant $\hat{A}_{q,r}$, in comparison with the lower bound $\widetilde A_{q,r}$.}
\label{fig:u0-kdv}
\end{figure}
Each optimal neural network $\phi_{\theta^*}(x)$ can be matched with a particular mKdV breather, as shown in the left panel of Figure \ref{fig:u0-kdv}. Comparing the value of $\frac\alpha\beta$ in each breather in Figure \ref{fig:u0-kdv}, we infer that all pairs try to fit the same breather, in this case, with $\alpha = 4$. Table \ref{tab:Aqr} summarizes the obtained optimal constant $\hat{A}_{q,r}$, averaged over 5 independent runs of the algorithm in each admissible pair. Additionally, the last column shows $\widetilde{A}_{q,r}$ for each admissible pair. As can be seen in the table, the values of $\hat{A}_{q,r}$ tend to stabilize around the lower bound $\widetilde{A}_{q,r}$, but even after $50\,000$ iterations they remain strictly below this lower bound.  Finally, Table \ref{tab:error_Aqr} shows the relative error of $\hat{A}_{q,r}$ with respect to $\widetilde{A}_{q,r}$. In each case, the relative error after $50\,000$ iterations is lower than $4 \times 10^{-3}$.

\begin{table}[!ht]
\centering
\caption{Evolution of $\hat A_{q,r}$ (iterations), compared with the Frank--Sabin lower bound $\widetilde A_{q,r}$. In every admissible pair $\hat A_{q,r}$ approaches $\widetilde A_{q,r}$ monotonically from below but never crosses it.}
\label{tab:Aqr}
\begin{tabular}{lccccc}
\toprule
Iterations & 100 & 1000 & 10000 & 50000 & $\widetilde A_{q,r}$ \\
\midrule
$\hat A_{5,10}$ & 0.5568 & 0.7779 & 0.7890 & 0.7915 & 0.7926 \\
$\hat A_{6,6}$  & 0.6168 & 0.7773 & 0.7852 & 0.7865 & 0.7886 \\
$\hat A_{8,4}$  & 0.6408 & 0.8012 & 0.8067 & 0.8084 & 0.8112 \\
$\hat A_{12,3}$ & 0.7102 & 0.8463 & 0.8535 & 0.8541 & 0.8556 \\
\bottomrule
\end{tabular}
\end{table}

\begin{table}[t]
\centering
\caption{Evolution of ${\bf error}_{q,r}$ (iterations), with target the Frank--Sabin lower bound $\widetilde A_{q,r}$.}
\label{tab:error_Aqr}
\begin{tabular}{lcccc}
\toprule
Iterations & 100 & 1000 & 10000 & 50000 \\
\midrule
${\bf error}_{5,10}$ & 2.975 $\times 10^{-1}$ & 1.858 $\times 10^{-2}$ & 4.546 $\times 10^{-3}$ & 1.642 $\times 10^{-3}$ \\
${\bf error}_{6,6}$  & 2.179 $\times 10^{-1}$ & 1.436 $\times 10^{-2}$ & 4.345 $\times 10^{-3}$ & 2.576 $\times 10^{-3}$ \\
${\bf error}_{8,4}$  & 2.101 $\times 10^{-1}$ & 1.228 $\times 10^{-2}$ & 5.596 $\times 10^{-3}$ & 3.385 $\times 10^{-3}$\\
${\bf error}_{12,3}$ & 1.699 $\times 10^{-1}$ & 1.055 $\times 10^{-2}$ & 2.398 $\times 10^{-3}$ & 1.680 $\times 10^{-3}$ \\
\bottomrule
\end{tabular}
\end{table}

\subsection{Airy-Strichartz for a non critical case} \label{app:AS-noncrit}
In this case we consider the pair $(q,r)=(8,8)$. For such a pair we have $\gamma = 0$. Extremizers for the Airy-Strichartz inequality in this case are well known to exist \cite{HundertmarkShao10}, but the extremizer profile is still unknown. We perform the AI-based algorithm to estimate the constant $\hat{A}_{8,8}$ and the extremizer profile. The results are summarized in Figure~\ref{fig:A88-result} and Table~\ref{tab:A88} for 10 independent seeds. From the figure one can see that all simulations yield to the same profile and the table shows that the sharp estimated constant is $0.7660$ with a negligible standard deviation. An $L^2$ function with the shape of the approximated profile has the form
\[
u_0(x) = A e^{ax^2} + B x^2e^{bx^2} + Cx^4e^{cx^2},
\]
with parameters $A = 1.0097$, $a = -0.1749$, $B = -1.3519$, $b = -0.6799$, $C = -0.0349$, $c = -0.3111$. A comparison of such a $u_0$ with the approximated profile can be found in Figure~\ref{fig:A88-fit}.
\begin{figure}[t]
\centering
\begin{subfigure}[t]{0.45\linewidth}\centering
\includegraphics[width=\linewidth]{images/u0-880.png}
\caption{Approximate extremizer $\phi_{\theta^\star}(x)$ for the pair $(8,8)$.}
\end{subfigure}\hfill
\begin{subfigure}[t]{0.45\linewidth}\centering
\includegraphics[width=\linewidth]{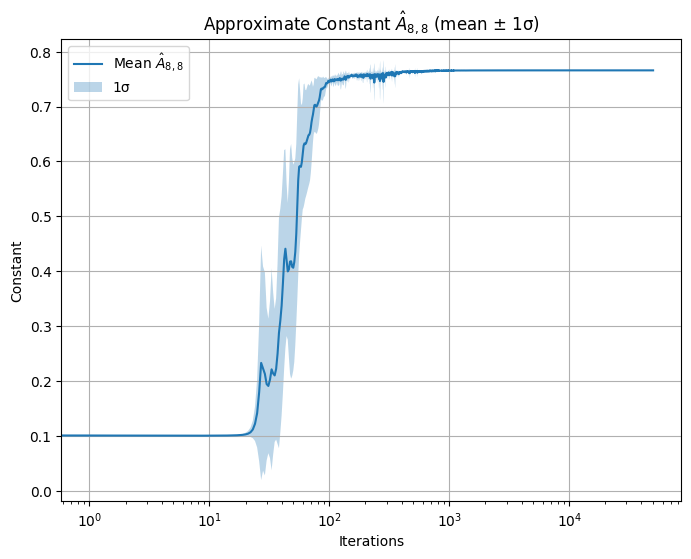}
\caption{Convergence of $\hat A_{8,8}$.}
\end{subfigure}
\caption{AI-based discovery for the Airy--Strichartz inequality, pair $(8,8)$.}
\label{fig:A88-result}
\end{figure}

\begin{table}[!ht]
\centering
\caption{Evolution of $\hat A_{8,8}$ (mean $\pm$ std) across iterations.}
\label{tab:A88}
\begin{tabular}{lc}
\toprule
Iterations & $\hat A_{8,8}$ \\
\midrule
100   & $0.7460 \pm 3.87 \times 10^{-3}$  \\
1000    & $0.7656 \pm 1.44 \times 10^{-4}$  \\
10000   & $0.7659 \pm 1.32 \times 10^{-6}$  \\
50000   & $0.7660 \pm 1.12 \times 10^{-5}$  \\
\bottomrule
\end{tabular}
\end{table}

\begin{figure}
    \centering
    \includegraphics[width=0.5\linewidth]{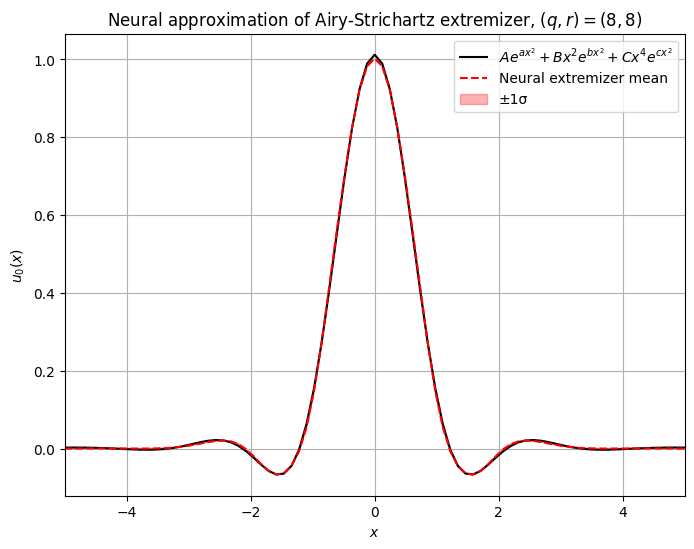}
    \caption{Best function fitting the approximated profile, for the pair $(8,8)$.}
    \label{fig:A88-fit}
\end{figure}

\subsection{Stability analysis}
For this purpose, we fix the pair to be $(q,r) = (6,6)$, and the other pairs behave in a similar way. We perform 10 independent realizations of the algorithm with parameters $R = 500$, $T = 5$, $N = 2^{13}$, $M = 2^{10}$ and $5\times 10^5$ iterations. Figures \ref{fig:A66_mean} and \ref{fig:A66_relative_error} summarize the mean $\pm$ the standard deviation. As can be seen in the figures, across all iterations the value of $\hat{A}_{6,6}$ does not reach the value of $\widetilde{A}_{6,6}$, and the standard deviation of the relative error is less than $5\times 10^{-3}$ after $5000$ iterations.

\begin{figure}
    \centering
    \includegraphics[width=0.5\linewidth]{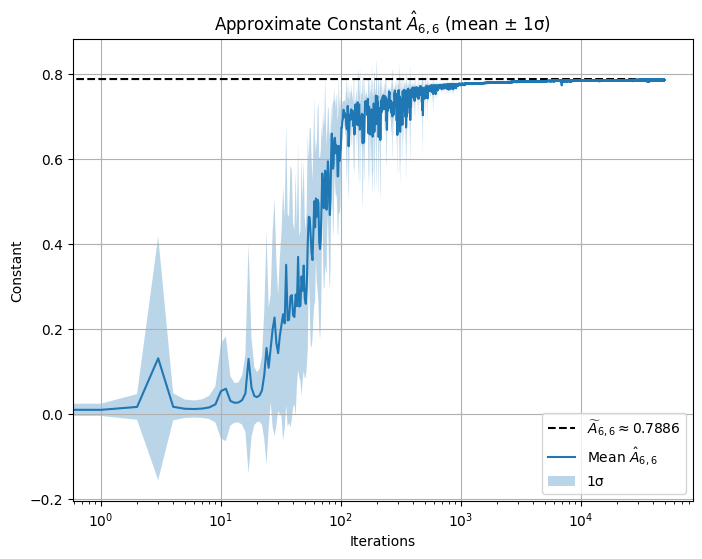}
    \caption{In blue, Estimation of the sharp constant with its standard deviation. In black, the value of $\widetilde{A}_{6,6}$.}
    \label{fig:A66_mean}
\end{figure}

\begin{figure}
    \centering
    \includegraphics[width=0.5\linewidth]{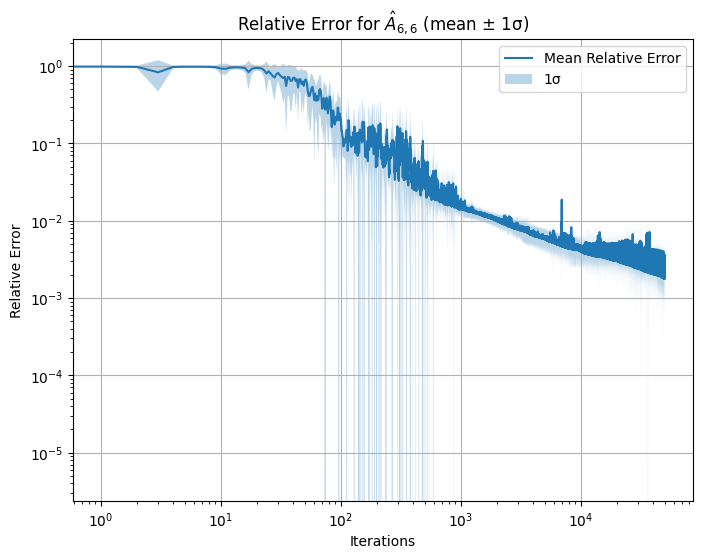}
    \caption{Relative error of $\hat A_{6,6}$ with respect to $\widetilde{A}_{6,6}$.}
    \label{fig:A66_relative_error}
\end{figure}

\subsection{Activation function analysis}
As an additional check we re-ran the same pipeline with $\tanh$ activation in place of the wavelet activation \eqref{act}, keeping all other hyperparameters fixed. Figure~\ref{fig:u0-kdv-tanh} shows the resulting profile and relative error over $10$ independent runs at $(q,r)=(6,6)$. The qualitative outcome is the same: the discovered profile is again an mKdV breather, the ratio approaches $\widetilde A_{6,6}$ from below, and the relative error stabilises at the $10^{-3}$ level. Quantitatively, the breather fitted with $\tanh$ has fewer internal oscillations (smaller $\alpha$) than with the wavelet activation. This is consistent with the picture that the absolute value of $\alpha$ at any finite training time reflects the spectral expressivity of the architecture, while the limit object (the breather family $\alpha\to\infty$) and the limit constant ($\widetilde A_{q,r}$) do not. The wavelet activation reaches higher effective $\alpha$ for the same compute budget, which is why we use it as default.
\begin{figure}[t]
\centering
\begin{subfigure}[t]{0.45\linewidth}\centering
\includegraphics[width=\linewidth]{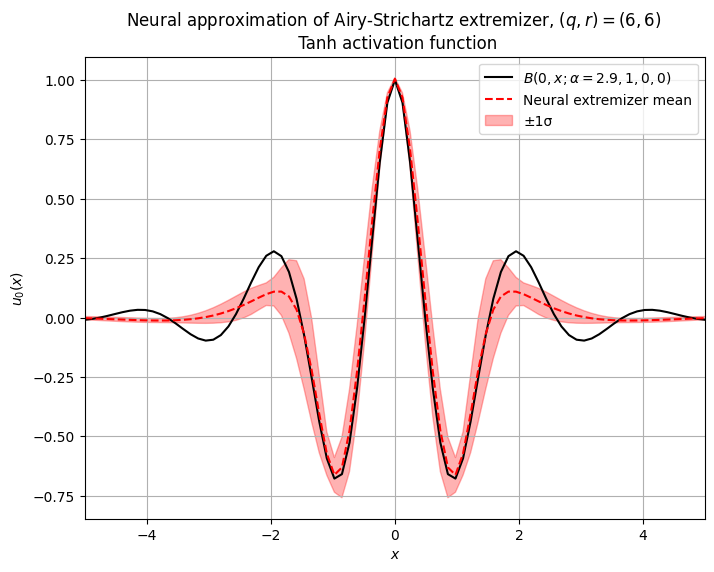}
\caption{Approximate extremizer $\phi_{\theta^\star}(x)$ for the pair $(6,6)$ with Tanh activation function.}
\end{subfigure}\hfill
\begin{subfigure}[t]{0.45\linewidth}\centering
\includegraphics[width=\linewidth]{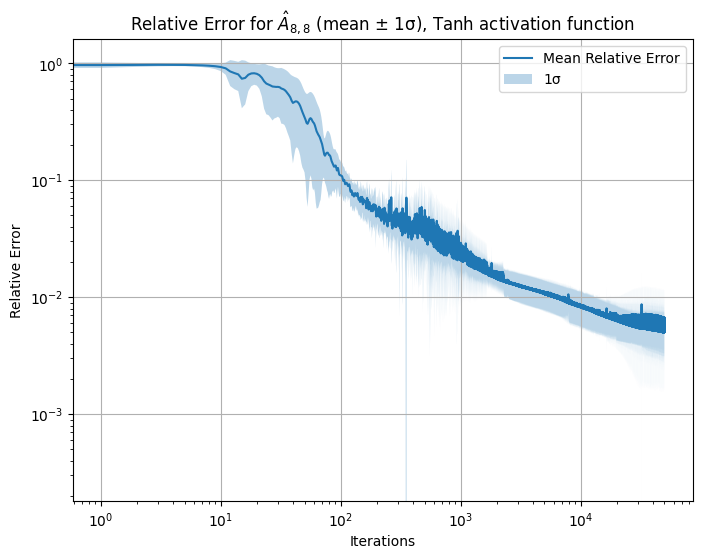}
\caption{Relative error of $\hat A_{6,6}$.}
\end{subfigure}
\caption{AI-based discovery for the Airy--Strichartz inequality, pair $(6,6)$.}
\label{fig:u0-kdv-tanh}
\end{figure}

The numerical experiments on the Airy--Strichartz inequalities below allow us to formulate the following conjecture, equivalent to the version stated in Conjecture~\ref{conj:main} of Section~\ref{sec:airy}.
\begin{conjecture}\label{conj:appendix}
 For every Schr\"odinger-admissible pair $(q,r)$ in the limiting case $\gamma=\frac1q$, the sharp Airy--Strichartz constant equals the Frank--Sabin lower bound,
\[
A_{q,r} \;=\; \widetilde A_{q,r} \;=\; 3^{-1/q}\,a_r\,S_{1,q,r},
\]
and is not attained by any $u_0\in L^2(\R)$. Moreover,
\[
A_{q,r} \;=\; \lim_{\alpha\to\infty}\, R[B(0,\cdot;\alpha,1,0,0)],
\]
modulo the symmetries of the Airy equation. In particular, the mKdV breather family $B(0,\cdot;\alpha,1,0,0)$ is an increasing (in $\alpha$) maximizing sequence approaching the Frank--Sabin bound from below; isolated breathers also appear as critical points of the ratio below $\widetilde A_{q,r}$ in several admissible pairs.
\end{conjecture}

\section{Neural network hyperparameters}
\label{app:hyper}

All experiments are implemented in PyTorch. Most runs were performed on a 64-bit MacBook Pro with an M4 chip (MPS backend, single precision). When the exponents in the mixed norm are large, single precision becomes unstable (overflow leading to NaN); in those cases we switched to double precision on Google Colab Pro with an NVIDIA T4 GPU. To quantify the precision effect, we ran five paired experiments at $(d,q,r)=(1,6,6)$: in float32 we obtained $\hat S_{1,6,6}\approx 0.8098$, while in float64 we obtained $\hat S_{1,6,6}\approx 0.8107$, an absolute difference of $9.17\times 10^{-4}$.

\paragraph{Airy experiments.} The spatial grid is $[-R,R]$ with $R=500$ and $N=2^{13}$ points, and the time grid is $[-T,T]$ with $T=5$ and $M=1024$ points. The network has $4$ hidden layers of $20$ neurons each and uses the wavelet activation \eqref{act} with the initialization of Section~\ref{sec:method}. Training uses Adam with a cyclic learning rate between $10^{-3}$ and $10^{-5}$, for up to $5\times 10^4$ iterations. The ridge is $\lambda=10^{-6}$, the $L^2$ normalization is applied at every forward pass, and the boundary term $|\phi_\theta(0)-1|^2$ is unweighted.

\paragraph{Schr\"odinger experiments.} The spatial grid is $[-R,R]$ with $R=30$, and $M=N=1024$ points (in $d=1$). In dimension $d=2$, $R=10$, $M=32$, $N=128$. Network as above ($4$ layers, $20$ neurons, wavelet activation). Adam with cyclic schedule, up to $10^4$ iterations. The ridge is $\lambda=10^{-6}$, the $L^2$ normalization is applied at every forward pass, and the boundary term $||\phi_\theta(0)|-1|^2$ is unweighted.

\paragraph{Random seeds and reproducibility.} For every experiment with reported error bars, we run $5$ (Schr\"odinger non-endpoint, Airy main) or $10$ (Airy stability, $(8,8)$ non-critical case) independent realizations differing only in the random seed of the network initialization.

\section{Outlook: analytic strategies and other dispersive groups}
\label{app:outlook}

\paragraph{A possible route to a proof of Conjecture~\ref{conj:main}.} A rigorous evaluation of $\lim_{\alpha\to\infty}R[B(0,\cdot;\alpha,1,0,0)]$ along the Airy flow, combined with the sharp Schr\"odinger constant $S_{1,q,r}$ of \cite{Foschi07,HundertmarkZharnitsky06}, should in principle saturate the Frank--Sabin bound \eqref{FS-bound}. The heat-flow monotonicity of \cite{BBC H09} and the integral identities of \cite{HundertmarkZharnitsky06} are plausible ingredients. The key technical step is an asymptotic analysis of the Strichartz integral $\|\,|D_x|^\gamma e^{-t\partial_x^3}B(0,\cdot;\alpha,1,0,0)\|_{L^q_tL^r_x}$ as $\alpha\to\infty$, which we intend to pursue in subsequent work.

\paragraph{Beyond Airy.} The neural optimization pipeline of Section~\ref{sec:method} is expected to apply without modification, up to the choice of propagator, to
\begin{itemize}
\item the wave propagator $e^{\pm it|\nabla|}$ on $\R^d$, $d\ge 2$;
\item the Klein--Gordon propagator $e^{it\sqrt{m^2-\Delta}}$;
\item the Zakharov--Kuznetsov propagator $U(t)=e^{t\partial_{x_1}\Delta}$ in $d\in\{2,3\}$.
\end{itemize}
In each case the analogous Frank--Sabin-type dichotomy and the identification of candidate maximizing sequences are open. Preliminary runs on the wave model in $d=2$ indicate the appearance of localized-oscillatory structures that are Gaussian-like; a systematic treatment is deferred to a forthcoming paper.



\end{document}